\documentclass{article}
\usepackage{rims02e}
\usepackage{amssymb,latexsym,amsfonts}

\title{Algebraic transformations of Gauss
hypergeometric functions}

\author{Raimundas Vid\=unas\thanks{Supported by the Dutch NWO project 613-06-565,
and by the Kyushu University 21 Century COE Programme "Development of Dynamic Mathematics
with High Functionality" of the Ministry of Education, Culture, Sports, Science
and Technology of Japan.}%
\affil{Faculty of Mathematics, Kobe University}%
\mail{vidunas@math.kobe-u.ac.jp}}

\newtheorem{theorem}{Theorem}[section]
\newtheorem{lemma}[theorem]{Lemma}

\newtheorem{remark}[theorem]{Remark}
\usepackage{latexsym}

\newcommand{\hpg}[5]{{}_{#1}\mbox{\rm F}_{\!#2}\!
  \left(\left.{#3 \atop #4}\right| #5 \right) }

\newcommand{\hpgo}[2]{{}_{#1}\mbox{\rm F}_{\!#2}}

\newcommand{\app}[4]{F_{\!#1}\!
  \left(\left.{#2 \atop #3}\right| #4 \right) }
\newcommand{\proof}{{\bf Proof. }}
\newcommand{\equal}{\!\!\!=\!\!\!}
\newcommand{\singset}{\Delta}
\newcommand{\CC}{{\Bbb C}}

\newcommand{\PP}{{\Bbb P}}
\newcommand{\RR}{{\Bbb R}}

\newcommand{\ZZ}{{\Bbb Z}}

\begin{document}

\maketitle

\begin{abstract} This article gives a classification scheme of algebraic transformations
of Gauss hypergeometric functions, or pull-back transformations between
hypergeometric differential equations. The classification recovers the
classical transformations of degree 2, 3, 4, 6, and finds other
transformations of some special classes of the Gauss hypergeometric
function. The other transformations are considered more thoroughly in 
a series of supplementing articles.
\end{abstract}

\section{Introduction}

An algebraic transformation of Gauss hypergeometric functions is an identity
of the form
\begin{equation} \label{hpgtransf}
\hpg{2}{1}{\!\widetilde{A},\,\widetilde{B}\,}{\widetilde{C}}{\,x}
=\theta(x)\;\hpg{2}{1}{\!A,\,B\,}{C}{\varphi(x)}.
\end{equation}
Here $\varphi(x)$ is a rational function of $x$, and $\theta(x)$ is a
radical function, i.e., a product of some powers of rational functions.
Examples of algebraic transformations are the following well-known quadratic
transformations (see \cite[Section~2.11]{bateman}, \cite[formulas 38,
45]{goursat}):
\begin{eqnarray} \label{quadr1}
\hpg{2}{1}{a,\,b}{\frac{a+b+1}{2}}{x} & = &
\hpg{2}{1}{\frac{a}{2},\,\frac{b}{2}}{\frac{a+b+1}{2}}{4x\,(1-x)},\\
\label{quadr2} \hpg{2}{1}{a,\;b}{2b}{x} & = &
\left(1-\frac{x}{2}\right)^{-a} \hpg{2}{1}{\frac{a}{2},\,\frac{a+1}{2}}
{b+\frac{1}{2}}{\frac{x^2}{(2-x)^2}}.
\end{eqnarray}

Algebraic transformations of Gauss hypergeometric functions are usually
induced by pull-back transformations between their hypergeometric
differential equations. General relation between these two kinds of
transformations is given in Lemma \ref{transeqv} here below. By that lemma,
if a pull-back transformation converts a hypergeometric equation to a
hypergeometric equation as well, then there are identities of the form
(\ref{hpgtransf}) between hypergeometric solutions of the two hypergeometric
equations. %, unless the transformed equation has a trivial monodromy group.
Conversely, an algebraic transformation (\ref{hpgtransf}) is induced by a
pull-back transformation of the corresponding hypergeometric equations,
unless the hypergeometric series on the left-hand side of (\ref{hpgtransf})
satisfies a first order linear differential equation.

This article classifies pull-back transformations between hypergeometric
differential equations. At the same time we essentially classify algebraic
transformations (\ref{hpgtransf}) of Gauss hypergeometric functions.
Classical fractional-linear and quadratic transformations are due to Euler,
Pfaff, Gauss and Kummer. In \cite{goursat} Goursat gave a list of
transformations of degree 3, 4 and 6. It has been widely assumed that there
are no other algebraic transformations, unless hypergeometric functions are
algebraic functions. For example, \cite[Section 2.1.5]{bateman} states the
following: ``Transformations of [degrees other than 2, 3, 4, 6] can exist
only if $a,b,c$ are certain rational numbers; in those cases the solutions
of the hypergeometric equation are algebraic functions." As our study shows,
this assertion is unfortunately not true. This fact is noticed in \cite{andkitaev} as well. 
Existence of a few special transformations follows from \cite{hodgkins1}, \cite{beukers}.

Regarding transformations of algebraic hypergeometric functions 
(or more exactly, pull-back transformations of hypergeometric differential 
equations with a finite monodromy group), celebrated Klein's theorem \cite{klein77}
ensures that all these hypergeometric equations are pull-backs of a few standard hypergeometric
equations. Klein's pull-back transformations do not change the {\em projective monodromy group}.
The possible finite projective monodromy groups are: a cyclic (including the trivial), 
a finite dihedral, the tetrahedral,  the octahedral or the icosahedral groups. 
Transformations of algebraic hypergeometric functions that reduce the projective monodromy group
are compositions of a few ``reducing" transformations and Klein's transformation
keeping the smaller monodromy group; see Remark \ref{trivexcept} below.

The ultimate list of pull-back transformations between hypergeometric
differential equations (and of algebraic transformations for their
hypergeometric solutions) is the following:
\begin{itemize}
\item Classical algebraic transformations of degree 2, 3, 4 and 6 due to
Gauss, Euler, Kummer, Pfaff and Goursat. %, Klein, Riemann, Goursat.
%These include fractional-linear transformations, quadratic transformations,
%and Goursat's transformations of degree 3, 4 and 6.
We review classical transformations in Section \ref{clalgtr},
including fractional-linear transformations.%
\item Transformations of hypergeometric equations with an abelian monodromy
group. This is a degenerate case \cite{degeneratehpg}; the hypergeometric equations have 2
(rather than 3) actual singularities. We consider these transformations in
Section \ref{logarithms}.%
\item Transformations of hypergeometric equations with a dihedral monodromy
group. These transformations are considered here in Section \ref{dihedrals},
or more thoroughly in \cite[Sections 6 and 7]{tdihedral}.%
\item Transformations of hypergeometric equations with a finite monodromy group. 
The hypergeometric solutions are algebraic functions. 
Transformations of hypergeometric equations with finite cyclic or dihedral monodromy
groups can be included in the previous two cases.
Transformations of hypergeometric equations with the
tetrahedral, octahedral or icosahedral projective monodromy groups 
are considered here in Section \ref{alggtranfs}, or more thoroughly in \cite{talggaus}.%
\item Transformations of hypergeometric functions which are incomplete
elliptic integrals. These transformations correspond %cohere
to endomorphisms of certain elliptic curves. They are considered in Section \ref{ellints},
or more thoroughly in \cite{telliptici}.
\item Finitely many transformations of so-called {\em hyperbolic
hypergeometric functions.} Hypergeometric equations for these functions have
local exponent differences $1/k_1$, $1/k_2$, $1/k_3$, where $k_1,k_2,k_3$
are positive integers such that $1/k_1+1/k_2+1/k_3<1$. These transformations
are described in Section \ref{otherat}, or more thoroughly
in \cite{thyperbolic}. % and in \cite{thyperbolic}.
\end{itemize}
The classification scheme is presented in Section \ref{clscheme}. Sections
\ref{clalgtr} through \ref{otherat} characterize various cases of algebraic
transformations of hypergeometric functions. We mention some three-term
identities with Gauss hypergeometric functions as well. The non-classical
cases are considered %(or will be considered)
more thoroughly in separate articles \cite{tdihedral}, \cite{talggaus}, 
\cite{telliptici}, \cite{thyperbolic}.

Recently, Kato \cite{kato} classified algebraic transformations of the $\hpgo32$ 
hypergeometric series. The rational transformations for the argument $z$ in that list
form a strict subset of the transformations considered here.

\section{Preliminaries}
\label{preliminar}

The hypergeometric differential equation is \cite[Formula
(2.3.5)]{specfaar}:
\begin{equation} \label{hpgde}
z\,(1-z)\,\frac{d^2y(z)}{dz^2}+
\big(C-(A\!+\!B\!+\!1)\,z\big)\,\frac{dy(z)}{dz}-A\,B\,y(z)=0.
\end{equation}
This is a Fuchsian equation %on the complex projective line $\PP^1$
with 3 regular singular points $z=0,1$ and $\infty$. % on $\PP^1$.
The local exponents are:
\[
\mbox{$0$, $1-C$ at $z=0$;} \qquad \mbox{$0$, $C-A-B$ at $z=1$;}\qquad
\mbox{and $A$, $B$ at $z=\infty$.}
\]
A basis of solutions for general equation (\ref{hpgde}) is
\begin{equation} \label{basissol}
\hpg{2}{1}{A,\,B\,}{C}{\,z}, \qquad
z^{1-C}\;\hpg{2}{1}{\!1+A-C,\,1+B-C}{2-C}{\,z}.
\end{equation}
For basic theory of hypergeometric functions and Fuchsian equations %we refer to
see \cite{beukers}, \cite[Chapters 1, 2]{alexaw} or \cite[Chapters 4,
5]{nicoboek}. We use the approach of Riemann and Papperitz \cite[Sections
2.3, 3.9]{specfaar}.

A {\em (rational) pull-back transformation} of an ordinary linear
differential equation %the hypergeometric equation
has the form
\begin{equation} \label{algtransf}
z\longmapsto\varphi(x), \qquad y(z)\longmapsto
Y(x)=\theta(x)\,y(\varphi(x)),
\end{equation}
where $\varphi(x)$ and $\theta(x)$ have the same meaning as in formula
(\ref{hpgtransf}). Geometrically, this transformation {\em pull-backs} a
differential equation %the hypergeometric equation
on the projective line $\PP^1_z$ to a differential equation on the
projective line $\PP^1_x$, with respect to the finite covering
$\varphi:\PP^1_x\to\PP^1_z$ determined by the rational function
$\varphi(x)$. %We use the notations $\PP^1_x$, $\PP^1_z$ 
Here and throughout the paper, we let $\PP^1_x$, $\PP^1_z$ denote the
projective lines with rational parameters $x$, $z$ respectively.
A pull-back transformation of a Fuchsian equations gives a Fuchsian
equation again. In \cite{andkitaev} pull-back transformations %(\ref{algtransf})
are called {\em $RS$-transformations}.

We introduce the following definition: an {\em irrelevant singularity} for
an ordinary differential equation is a regular singularity which is not
logarithmic, and where the local exponent difference is equal to 1. 
An irrelevant singularity can
be turned into a non-singular point after a suitable pull-back
transformation (\ref{algtransf}) with $\varphi(x)=x$.
(For comparison,  an {\em apparent singularity} is a regular singularity which is
not logarithmic, and where the local exponents are integers. Recall that at
a {\em logarithmic point} is a singular point where there is only one local
solution of the form
$x^\lambda\!\left(1+\alpha_1x+\alpha_2x^2+\ldots\right)$, where $x$ is a
local parameter there.) For us, a {\em relevant singularity} is a singular
point which is not an irrelevant singularity. 

We are interested in pull-back transformations of one hypergeometric
equation to other hypergeometric equation, possibly with different
parameters $A,B,C$. These pull-back transformations are related to algebraic
transformations of Gauss hypergeometric functions as follows.
\begin{lemma} \label{transeqv}
\begin{enumerate}
\item Suppose that pull-back transformation $(\ref{algtransf})$ of
hypergeometric equation $(\ref{hpgde})$ is a hypergeometric equation as well
(with the new indeterminate $x$).
%, and that the transformed equation has non-trivial monodromy. 
Then, possibly after fractional-linear transformations on %the projective lines
$\PP_x^1$ and $\PP^1_z$, there is an identity of the form
$(\ref{hpgtransf})$ between hypergeometric solutions of the two
hypergeometric equations.%
\item Suppose that hypergeometric identity $(\ref{hpgtransf})$ holds in some
region of the complex plane. Let $Y(x)$ denote the left-hand side of the
identity. If %the logarithmic derivative
$Y'(x)/Y(x)$  is not a rational function of $x$, then the transformation
$(\ref{algtransf})$ converts the hypergeometric equation $(\ref{hpgde})$
into a hypergeometric equation for $Y(x)$.
\end{enumerate}
\end{lemma}
\proof We have a two-term identity whenever we have a singular point $S\in\PP_x^1$
of the transformed equation above a singular point $Q\in\PP_z^1$ of the starting equation.
Using fractional-linear transformations on $\PP_x^1$ and $\PP^1_z$
we can achieve $S$ is the point $x=0$ and that $Q$ is the point $z=0$. 
Then identification of two hypergeometric solutions with the local exponent 0 
and the value 1 at (respectively) $x=0$ and $z=0$ gives a two-term identity
as in (\ref{hpgtransf}). If all three singularities of the transformed equation
do not lie above $\{0,1,\infty\}\subset\PP_z^1$, they are apparent singularities.
Then the transformed equation has trivial monodromy, while the starting hypergeometric
equation has a finite monodromy group. As we will consider explicitly in Sections
\ref{alggtranfs} and \ref{logarithms}, \ref{dihedrals}, the pull-back transformations
reducing the monodromy group and the pull-back transformations keeping
the trivial monodromy group can be realized by two-term hypergeometric identities.
This is recaped in Remark \ref{trivexcept} below.

For the second statement, we have two second-order linear differential
equations for the left-hand side of (\ref{hpgtransf}): the hypergeometric
equation for $Y(x)$, and a pull-back transformation (\ref{algtransf}) of the
hypergeometric equation (\ref{hpgde}). If these two equations are not
$\CC(x)$-proportional, then we can combine them linearly to a first-order
differential equation $Y'(x)=r(x)\,Y(x)$ with $r(x)\in\CC(x)$,
contradicting %%%which would contradict
the condition on $Y'(x)/Y(x)$. \qed\\

\noindent
If we have an identity (\ref{hpgtransf}) without a pull-back transformation
between corresponding hypergeometric equations, the left-hand side of the
identity can be expressed as terminating hypergeometric series up to a
radical factor; see Kovacic algorithm \cite{kovacic}, \cite[Section
4.3.4]{vdputsing}. In a formal sense, any pair of terminating hypergeometric series
is algebraically related. We do not consider these degenerations.
\begin{remark} \rm
We also do not consider transformations of the type
$\hpgo{2}{1}(\varphi_1(z))=\theta(z)\,\hpgo{2}{1}(\varphi_2(z))$, where
$\varphi_1(z)$, $\varphi_2(z)$ are rational functions (of degree at least
2). Therefore we miss transformations of some complete elliptic integrals,
such as
\begin{equation} \label{elliptick}
{\bf K}(x)=\frac{2}{1+y}\;{\bf K}\!\left(\frac{1-y}{1+y}\right),
%\qquad \mbox{where } x^2+y^2=1,
\end{equation}
where $x^2+y^2=1$ and
\[ %begin{equation}
{\bf K}(x)=\frac{\pi}{2}\,\hpg{2}{1}{1/2,\,1/2}{1}{\,x^2}=
{\displaystyle\int_0^1} \frac{dt}{\sqrt{(1-t^2)(1-x^2t^2)}}.
\] %end{equation}
Identity (\ref{elliptick}) plays a key role in the theory of
arithmetic-geometric mean; see \cite[Chapter 3.2]{specfaar}. Other similar
example is the following formula, proved in \cite[Theorem 2.3]{bbgramanuj}:
\begin{equation}
\hpg{2}{1}{c,\,c+\frac{1}{3}\,}{\frac{3c+5}{6}}{\,x^3} = \big( 1+2x
\big)^{-3c}\; \hpg{2}{1}{c,\,c\!+\!\frac{1}{3}} {\frac{3c+1}{2}}
{1-\frac{(1\!-\!x)^3}{(1\!+\!2x)^3}}.
\end{equation}
The case $c=1/3$ was found earlier in \cite{borweins91}.
\end{remark}

A pull-back transformation between hypergeometric equations usually gives
several identities like (\ref{hpgtransf}) between some of the 24 Kummer's
solutions of both equations. It is appropriate to look first for
suitable pull-back coverings $\varphi:\PP_x^1\to\PP_z^1$ up to fractional-linear
transformations. As we will see, suitable pull-back coverings are determined
by appropriate transformations of singular points and local exponent differences.

Once a suitable covering $\varphi$ is known, it is
convenient to use Riemann's $P$-notation for deriving hypergeometric
identities (\ref{hpgtransf}) with the argument $\varphi(x)$. Recall that a
Fuchsian equation with 3 singular points is determined by the location of
those singular points and local exponents there. The linear space of
solutions is determined by the same data. It can be defined homologically without
reference to hypergeometric equations as a {\em local system} on the
projective line; see \cite{katzls}, \cite[Section 1.4]{graylde}. The
notation for it is
\begin{equation} \label{gensols}
P\left\{\begin{array}{ccc} \alpha & \beta & \gamma \\ a_1 & b_1 & c_1 \\
a_2 & b_2 & c_2 \end{array} \;z\; \right\},
\end{equation}
where $\alpha,\beta,\gamma\in\PP^1_z$ are the singular points, and
$a_1,a_2$; $b_1,b_2$; $c_1,c_2$ are the local exponents at them,
respectively. Recall that second order Fuchsian equations with 3 relevant singularities are defined
uniquely by their singularities and local exponents, unlike general Fuchsian
equations with more than 3 singular points.
Our approach can be entirely formulated in terms of local systems, 
without reference to hypergeometric equations and their pull-back transformations. 
By Papperitz' theorem
\cite[Theorem~2.3.1]{specfaar} we must have
\[ a_1+a_2+b_1+b_2+c_1+c_2=1.
\]
We are looking for transformations of local systems of the form
\begin{equation} \label{hpgsols}
P\left\{\begin{array}{ccc} 0 & 1 & \infty \\ 0 & 0 & \widetilde{A} \\
1-\widetilde{C} & \widetilde{C}-\widetilde{A}-\widetilde{B} & \widetilde{B}
\end{array} \;x\; \right\} = \theta(x)\;P\left\{\begin{array}{ccc} 0 & 1 & \infty \\ 0 & 0 & A \\
1-C & C-A-B & B \end{array} \;\varphi(x)\; \right\} .
\end{equation}
The factor $\theta(x)$ should shift local exponents at irrelevant
singularities to the values 0 and 1, and it should shift one local exponent
at both $x=0$ and $x=1$ to the value 0. In intermediate computations,
Fuchsian equations with more than 3 singular points naturally occur, but
those extra singularities are irrelevant singularities. We extend Riemann's
$P$-notation and write
\begin{equation} \label{extendedp}
P\left\{ \begin{array}{cccccc} \alpha & \beta & \gamma & S_1 & \ldots & S_k \\
a_1 & b_1 & c_1 & e_1 & \ldots & e_k \\
a_2 & b_2 & c_2 & e_1+1 & \ldots & e_k+1 \end{array} \;z\;\right\}
\end{equation}
to denote the local system (of solutions of a Fuchsian equation) with
irrelevant singularities $S_1,\ldots,S_k$. This notation makes sense if a
local system exists (i.e., if the local exponents sum up to the right
value); then it can be transformed to a local system like (\ref{gensols}).
For example, if none of the points $\gamma$, $S_1,\ldots,S_k$ is the infinity, local system
(\ref{extendedp}) can be identified with
\[
\frac{(z-S_1)^{e_1}\ldots(z-S_k)^{e_k}}{(z-\gamma)^{e_1+\ldots+e_k}}\;
P\left\{\begin{array}{ccc} \alpha & \beta & \gamma \\ a_1 & b_1 & c_1+e_1+\ldots+e_k \\
a_2 & b_2 & c_2+e_1+\ldots+e_k \end{array} \;z\; \right\}.
\]
Here is an example of computation with local systems leading %which leads 
to quadratic transformation (\ref{quadr2}):
\begin{eqnarray*}
P\left\{ \begin{array}{ccc} 0&1&\infty \\ 0 & 0 & \frac{a}{2} \\
\frac{1}{2}-b&b-a&\frac{a+1}{2} \end{array} \,t^2\, \right\}
& = & P\left\{ \begin{array}{cccc} 0 & 1 & -1 & \infty \\ 0 & 0 & 0 & a \\
1-2b & b-a & b-a & a+1 \end{array} \;t\;\right\} \\
& = & P\left\{ \begin{array}{cccc} 0 & 1 & \infty & 2 \\ 0 & 0 & 0 & a \\
1-2b & b-a & b-a & a+1 \end{array} \;x=\frac{2t}{t+1}\,\right\} \\
& = & \left(2-x\right)^a\; P\left\{ \begin{array}{ccc} 0 & 1 & \infty \\
0 & 0 & a \\ 1-2b & b-a & b \end{array} \;x\; \right\}.
\end{eqnarray*}
To conclude (\ref{quadr2}), one has to identify two functions with the local
exponent 0 and the value 1 at $t=0$ and $x=0$ (in the first and the last
local systems respectively), like in the proof of part 1 of Lemma
\ref{transeqv}.

Once a hypergeometric identity (\ref{hpgtransf}) is obtained, it can be composed with
Euler's and Pfaff's fractional-linear transformations; we recall them in formulas 
(\ref{flinear1})--(\ref{flinear3}) below. % \cite[Theorem 2.2.5]{specfaar}. 
Geometrically, these transformations permute the 
singularities $1$, $\infty$ (on $\PP^1_z$ or $\PP^1_x$) and their local exponents.
Besides, simultaneous permutation of the local exponents at $x=0$ and $z=0$
usually implies a similar identity to (\ref{hpgtransf}), as presented in the following lemma. 
%We present possible hypergeometric identities up to the transformations mentioned just now.
\begin{lemma} \label{basislem}
Suppose that a pull-back transformation induces identity $(\ref{hpgtransf})$
in an open neighborhood of $x=0$. Then $\varphi(x)^{1-C}\sim K
x^{1-\widetilde{C}}$ as $x\to 0$ for some constant $K$, and the following
identity holds (if both hypergeometric functions are well-defined):
\begin{equation} \label{altbasis}
\hpg{2}{1}{\!1+\widetilde{A}-\widetilde{C},\,1+\widetilde{B}-\widetilde{C}}
{2-\widetilde{C}}{\,x\,}=\theta(x)\,\frac{\varphi(x)^{1-C}}{K\,x^{1-\widetilde{C}}}
\;\hpg{2}{1}{\!1+A-C,\,1+B-C}{2-C}{\varphi(x)}.
\end{equation}
\end{lemma}
\proof The asymptotic relation $\varphi(x)^{1-C}\sim K x^{1-\widetilde{C}}$ as $x\to 0$ 
is clear from the transformation of local exponents. (We are ensured that $\theta(0)=1$.)
Further, we have relation
(\ref{hpgsols}) and the relations
\begin{eqnarray*}
P\left\{\begin{array}{ccc} 0 & 1 & \infty \\ 0 & 0 & A \\
1-C & C-A-B & B \end{array} \;\varphi(x)\; \right\} & = & \varphi(x)^{1-C}\,
P\left\{\begin{array}{ccc} 0 & 1 & \infty \\ C-1 & 0 & A+1-C \\
0 & C-A-B & B+1-C \end{array} \;\varphi(x)\; \right\}, \\
P\left\{\begin{array}{ccc} 0 & 1 & \infty \\ 0 & 0 & \widetilde{A} \\
1-\widetilde{C} & \widetilde{C}-\widetilde{A}-\widetilde{B} & \widetilde{B}
\end{array} \;x\; \right\} & = & x^{1-\widetilde{C}}\;
P\left\{\begin{array}{ccc} 0 & 1 & \infty \\ \widetilde{C}-1 & 0 & \widetilde{A}+1-\widetilde{C} \\
0 & \widetilde{C}-\widetilde{A}-\widetilde{B} &
\widetilde{B}+1-\widetilde{C} \end{array} \;x\; \right\}.
\end{eqnarray*}
From here we get the right identification of local systems for
(\ref{altbasis}).\qed\\

A general pull-back transformation converts a hypergeometric equation to a
Fuchsian differential equation with several singularities. To find proper candidates for
pull-back coverings $\varphi:\PP_x^1\to\PP_z^1$, we look first for possible pull-back transformations
of hypergeometric equations to Fuchsian equations with (at most) 3 relevant
singularities. These Fuchsian equations can be always transformed to
hypergeometric equations by suitable fractional-linear pull-back
transformations, and vice versa. Relevant singular points and local exponent
differences for the transformed equation are determined by the covering
$\varphi$ only. Here are simple rules which determine singularities and
local exponent differences for the transformed equation.
\begin{lemma} \label{genrami}
Let $\varphi:\PP^1_x\to\PP^1_z$ be a finite covering. Let $H_1$ denote a
Fuchsian equation on $\PP^1_z$, and let $H_2$ denote the pull-back
transformation of $H_1$ under $(\ref{algtransf})$. Let $S\in\PP^1_x$,
$Q\in\PP^1_z$ be points such that $\varphi(S)=Q$.
\begin{enumerate}
\item The point $S$ is a logarithmic point for $H_2$ if and only if the
point $Q$ is a logarithmic point for $H_1$.%
\item If the point $Q$ is non-singular for $H_1$, then the point $S$ is not
a relevant singularity for $H_2$ if and only if the covering $\varphi$ does
not branch at $S$.%
\item If the point $Q$ is a singular point for $H_1$, then the point $S$ is
not a relevant singularity for $H_2$ if and only if the following two
conditions hold:
\begin{itemize}
\item The point $Q$ is not logarithmic.%
\item The local exponent difference at $Q$ is equal to $1/k$, where $k$ is
the branching index of $\varphi$ at $S$.%
\end{itemize}
\end{enumerate}
\end{lemma}
\proof First we note that if the point $S$ is not a relevant singularity,
then it is either a non-singular point or an irrelevant singularity.
Therefore $S$ is not a relevant singularity if and only if it is not a
logarithmic point and the local exponent difference is equal to 1.

Let $p$, $q$ denote the local exponents for $H_1$ at the point $Q$. Let $k$
denote the branching order of $\varphi$ at $S$. Then the local exponent
difference at $S$ is equal to $k(p-q)$. To see this, note that if $m\in\CC$
is the order of $\theta(x)$ at $S$, the local exponents at $S$ are equal to
$kp+m$ and $kq+m$. This fact is clear if $Q$ is not logarithmic, when the
local exponents can be read from solutions. In general one has to use the
indicial polynomial to determine local exponents.

The first statement is clear, since local solutions of $H_1$ at $S$ can be
pull-backed to local solutions of $H_2$ at $Q$, and local solutions of $H_2$
at $Q$ can be push-forwarded to local solutions of $H_1$ at $S$.

If the point $Q$ is non-singular, the point $S$ is not logarithmic by the
first statement, so $S$ is a not a relevant singularity if and only if
$k=1$. %(If $k\neq 1$ then $S$ is a relevant but apparent singularity.)

If the point $Q$ is singular, then the local exponent difference at $S$ is
equal to 1 if and only if the local exponent difference $|p-q|$ is equal to
$1/k$. \qed\\

The following Lemma gives an estimate for the number of points $S$ to which
part 3 of Lemma \ref{genrami} applies, and it gives a relation between local
exponent differences of two hypergeometric equations related by a pull-back
transformation and the degree of the pull-back transformation. In this paper
we make the convention that real local exponent differences are non-negative, and complex
local exponent differences have the argument in the interval $(-\pi,\pi]$.
\begin{lemma} \label{genrami2}
Let $\varphi:\PP^1_x\to\PP^1_z$ be a finite covering, and let $d$ denote the
degree of $\varphi$.
\begin{enumerate}
\item Let $\singset$ denote a set of $3$ points on $\PP^1_z$. If all
branching points of $\varphi$ lie above $\singset$, then there are exactly
$d+2$ distinct points on $\PP^1_x$ above $\singset$. Otherwise there are
more than $d+2$ distinct points above
$\singset$.%
\item Let $H_1$ denote a hypergeometric equation on $\PP^1_z$, and let $H_2$
denote a pull-back transformation of $H_1$ with respect to $\varphi$.
Suppose that $H_2$ is hypergeometric equation as well. Let $e_1,e_2,e_3$
%$e_0,e_1,e_\infty$ ... at $z=0$, $z=1$ and $z=\infty$ respectively
denote the local exponent differences for $H_1$, and let $e'_1,e'_2,e'_3$
denote the local exponent differences for $H_2$. Then
\begin{equation} \label{trareas}
d\left( e_1+e_2+e_3-1 \right)=e'_1+e'_2+e'_3-1.
\end{equation}
\end{enumerate}
\end{lemma}
\proof For a point $S\in\PP_x^1$ let $\mbox{ord}_S\,\varphi$ denote the
branching order of $\varphi$ at $S$. By Hurwitz formula \cite[Corollary
IV.2.4]{harts} we have $-2=-2d+D$, where
\[
D=\sum_{S\in\PP_x^1} \left( \mbox{ord}_S\,\varphi - 1 \right).
\]
Therefore $D=2d-2$. The number of points above $\singset$ is %equal to
\[
3d-\!\sum_{\varphi(S)\in\singset} \left( \mbox{ord}_S\,\varphi - 1 \right)
\,\ge\, 3d-D\,=\,d+2.
\]
We have the equality if and only if all branching points of $\varphi$ lie
above $\singset$.

Now we show the second statement. For a point $S\in\PP_z^1$ or
$S\in\PP_x^1$, let $\mbox{led}(S)$ denote the local exponent difference for
$H_1$ or $H_2$ (respectively) at $S$. The following sums make sense:
\begin{eqnarray*}
\sum_{S\in\PP_x^1} \big( \mbox{led}(S)-1 \big) & = &
\sum_{Q\in\PP_z^1}\,\sum_{\varphi(S)=Q} \big( \mbox{led}(S)-1 \big)\\
&=&\sum_{Q\in\PP_z^1}\left( d\;\mbox{led}(Q)-\sum_{\varphi(S)=Q}1 \right)\\
&=&d\sum_{Q\in\PP_z^1}\big(\mbox{led}(Q)-1\big)\,+\,D.
%\sum_{P\in\PP_x^1}\left( \mbox{ord}_P\,\varphi - 1 \right).
\end{eqnarray*}
The first sum is equal to $e'_1+e'_2+e'_3-3$. The last expression is equal
to $d\,(e_1+e_2+e_3-3)+2d-2$. \qed

\section{The classification scheme} \label{clscheme}

The core problem is to classify pull-back transformations of hypergeometric
equations to Fuchsian equations with at most 3 relevant singular points. By
Lemma \ref{genrami}, a general pull-back transformation gives a Fuchsian
equation with quite many relevant singular points, especially above the set
$\{0,1,\infty\}\subset\PP^1_z$. In order to get a Fuchsian equation with so
few singular points, we have to restrict parameters (or local exponent
differences) of the original hypergeometric equation, and usually we can
allow branching only above the set $\{0,1,\infty\}\subset\PP^1_z$.

We classify pull-back transformations between hypergeometric equations (and
algebraic transformations of Gauss hypergeometric functions) in the
following five principal steps:
\begin{enumerate}
\item Let $H_1$ denote hypergeometric equation (\ref{hpgde}), and consider
its pull-back transformation (\ref{algtransf}). Let $H_2$ denote the
pull-backed differential equation, and let $T$ denote the number of singular
points of $H_2$. Let $\singset$ denote the subset $\{0,1,\infty\}$ of
$\PP^1_z$,  and let $d$ denote the degree of the covering
$\varphi:\PP^1_x\to\PP_z^1$ in transformation (\ref{algtransf}). We
consequently assume that exactly $N\in\{0,1,2,3\}$ of the 3 local exponent
differences for $H_1$ at $\singset$ are {\em restricted} to the values
of the form $1/k$, where $k$ is a positive integer. If $k=1$ then 
the corresponding point of $\singset$ is assumed to be not logarithmic,
as we cannot get rid of singularities above a logarithmic point.
\item  In each assumed case, use Lemma \ref{genrami} and determine all %{\em a priori}
possible combinations of the degree $d$ and restricted local exponent differences. 
Let $k_1,\ldots,k_N$ denote the denominators of the restricted differences. 
By part 4 of Lemma \ref{genrami},
\begin{eqnarray*}
T&\ge& \mbox{[the number of singular points above $\singset$]}\\
&\ge& d+2-\mbox{[the number of non-singular points above $\singset$]}\\
&\ge& d+2-\sum_{j=1}^N \left\lfloor \frac{d}{k_j} \right\rfloor.
\end{eqnarray*}
Since we wish $T\le 3$, we get the following restrictive inequality in
integers:
\begin{equation} \label{nsingpoints}
d-\sum_{j=1}^N \left\lfloor \frac{d}{k_j} \right\rfloor\le 1.
\end{equation}
To skip specializations of cases with smaller $N$, we may assume that
$d\ge\max(k_1,\ldots,k_N)$. A preliminary list of possibilities can be
obtained by dropping the rounding down in (\ref{nsingpoints}); this gives a
weaker but more convenient inequality
\begin{equation} \label{nsingpoints2}
\frac1d\,+\,\sum_{j=1}^N \frac{1}{k_j}\,\ge\, 1.
\end{equation}
\item For each combination of $d$ and restricted local exponent differences,
determine possible branching patterns for $\varphi$ such that the
transformed equation $H_2$ would have at most three singular points. In most
cases we can allow branching points only above $\singset$, and
we have to take the maximal number $\left\lfloor d/k_j %\frac{d}{k_j}
\right\rfloor$ of non-singular points above the point with the local
exponent difference $1/k_j$.%
\item For each possible branching pattern, determine all rational functions
$\varphi(x)$ which define a covering with that branching pattern. For $d\le
6$ this can be done using a computer by a straightforward %naive 
method of undetermined coefficients. In \cite[Section 3]{thyperbolic} a more appropriate algorithm
is introduced which uses differentiation of $\varphi(x)$. In many cases this
problem has precisely one solution up to fractional-linear transformations.
But not for any branching pattern a covering exists, and there can be
several different coverings with the same branching pattern. For infinite
families of branching patterns we are able to give a general, algorithmic or
explicit characterization of corresponding coverings. For instance, if
hypergeometric solutions can be expressed very explicitly, we can %usually
identify the local systems in (\ref{hpgsols}) up to unknown factor
$\theta(x)$. Then quotients of corresponding hypergeometric solutions 
(aka {\em Schwarz maps}) can be identified precisely, which gives a straightforward way
to determine $\varphi(x)$. 
\item Once a suitable covering $\varphi:\PP^1_x\to\PP_z^1$ is computed,
there always exist corresponding pull-back transformations. Two-term
identities like (\ref{hpgtransf}) can be computed using extended Riemann's
$P$-notation of Section \ref{preliminar}.  We have two-term identities for 
each singular point $S$ of the transformed equation such that \mbox{$\varphi(S)\subset\Delta$},
as in the proof of part 1 of Lemma \ref{transeqv}. Once we fix $S$, $\varphi(S)$ as $x=0$,
$z=0$ respectively, permutations of local exponents and other singularities
give identities (\ref{hpgtransf}) which are related by Euler's and Pfaff's
transformations and Lemma \ref{basislem}. If the transformed equation has
less than 3 actual singularities, one can consider any point above $\Delta$
in this manner. Some of the obtained identities may be the same up to change
of free parameters.
\end{enumerate}
Now we sketch explicit appliance of the above procedure.
When $N=0$, i.e., when no local exponent differences are restricted, then
$d=1$ by formula (\ref{nsingpoints2}). This gives Euler's and
Pfaff's fractional-linear transformations. % \cite[Theorem 2.2.5]{specfaar}.
When $N=1$, we have the following cases:
\begin{itemize}
\item $k_1=2$, $d=2$. This gives the classical quadratic transformations.
See Section \ref{clalgtr}. %\cite[Section 3.9]{specfaar}.%
\item $k_1=1$, $d$ any. The $z$-point with the local exponent difference
$1/k_1$ is assumed to be non-logarithmic, so the equation $H_1$ has only two
relevant singularities. As we show in Lemma \ref{logpoint1} below,
the two unrestricted local exponent differences must be equal.
As it turns out, the covering $\varphi$ branches only above the two points
with unrestricted local exponent differences. If the triple of local
exponent differences for $H_1$ is $(1,p,p)$, the triple of local exponent
differences for $H_2$ is $(1,dp,dp)$. Formally, this case has a continuous
family of fractional-linear pull-back transformations, but that does not
give interesting hypergeometric identities.
\end{itemize}
When $N=2$, we have the following cases:
\begin{itemize}
\item If $\max(k_1,k_2)>2$, the possibilities are listed in Table
\ref{clasfig}. Steps 2 and 3 of the classification scheme are straightforward,
and a snapshot of possibilities after them is presented by the first four columns
of Table \ref{clasfig}. The notation for a branching pattern in the fourth column gives $d+2$
branching orders for the points above $\singset$; branching orders at points
in the same fiber are separated by the + signs, branching orders for
different fibers are separated by the = signs. Step 4 of our scheme gives at
most one covering (up to fractional-linear transformations) for each
branching pattern. Ultimately, Table \ref{clasfig} yields precisely the
classical transformations of degree 3, 4, 6 due to Goursat \cite{goursat};
see Section \ref{clalgtr}. It is straightforward to figure out possible
compositions of small degree coverings, and then identify them with the
unique coverings for Table \ref{clasfig}. Degrees of constituents for
decomposable coverings are listed in the last column from right (for the
constituent transformation from $H_1$) to left. Note that one degree 6
covering has two distinct decompositions; a corresponding hypergeometric
transformation is given in formula (\ref{eq:trd6}) below.
\begin{table} \begin{center}
\begin{tabular}{|c|c|c|c|c|}
\hline \multicolumn{2}{|c|}{Local exponent differences} & Degree & Branching
pattern above & Covering \\ \cline{1-2}
$(1/k_1,\,1/k_2,\,p)$ & above & $d$ & the regular singular points & composition \\
 \hline
$(1/2,\,1/3,\,p)$ & $(1/2,\,p,\,2p)$ & 3 & $2+1=3=2+1$ & indecomposable \\
$(1/2,\,1/3,\,p)$ & $(1/3,\,p,\,3p)$ & 4 & $2+2=3+1=3+1$ & indecomposable \\
$(1/2,\,1/3,\,p)$ & $\!(1/3,\,2p,\,2p)\!$ & 4 & $2+2=3+1=2+2$ & no covering \\
$(1/2,\,1/3,\,p)$ & $(p,\,p,\,4p)$ & 6 &
$2+2+2=3+3=4+1+1$ & $2\times 3$ \\
$(1/2,\,1/3,\,p)$ & $(2p,\,2p,\,2p)$ & 6 &
$2+2+2=3+3=2+2+2$ & $2\times3$ or $3\times2$ \\
$(1/2,\,1/3,\,p)$ & $(p,\,2p,\,3p)$ & 6 &
$2+2+2=3+3=3+2+1$ & no covering \\
$(1/2,\,1/4,\,p)$ & $(p,\,p,\,2p)$ & 4 & $2+2=4=2+1+1$ & $2\times 2$ \\
$(1/3,\,1/3,\,p)$ & $(p,\,p,\,p)$ & 3 & $3=3=1+1+1$ & indecomposable \\
\hline
\end{tabular}
\caption{Transformations of hypergeometric functions with 1 free parameter}
\label{clasfig} \end{center}\vspace{-2pt}
\end{table}
\item $k_1=2$, $k_2=2$, $d$ any. The monodromy group of $H_1$ is a dihedral
group. The hypergeometric functions can be expressed very explicitly, see
Section \ref{dihedrals}. The triple $(1/2,1/2,p)$ of local exponent
differences for $H_1$ is transformed either
to $(1/2,1/2,dp)$ for any $d$, or to $(1,dp/2,dp/2)$ for even $d$.%
\item $k_1=1$; $k_2$ and $d$ are any positive integers. The $z$-point with
the local exponent difference $1/k_1$ is not logarithmic, so the triple of
local exponent differences for $H_1$ must be $(1,1/k_2,1/k_2)$. The
monodromy group is a finite cyclic group. Possible transformations are
outlined in Section \ref{logarithms}.
\end{itemize}
When $N=3$, we have the following three very distinct cases:
\begin{itemize}
\item $1/k_1+1/k_2+1/k_3>1$. The monodromy groups of $H_1$ and $H_2$ are
finite, the hypergeometric functions are algebraic. The degree $d$ is
unbounded. Klein's theorem \cite{klein77} implies that any hypergeometric
equation with a finite monodromy group (or equivalently, with algebraic
solutions) is a pull-back transformation of a {\em standard hypergeometric
equation} with the same monodromy group. These are the most interesting
pull-back transformations for this case. Equations with finite cyclic
monodromy groups are mentioned in the previous subcase; their
transformations are considered in Section \ref{logarithms}. Equations with
finite dihedral monodromy groups are considered in Section \ref{dihedrals}.
Equations with the tetrahedral, octahedral or icosahedral projective 
monodromy groups are characterized in Section \ref{alggtranfs}.%
\item $1/k_1+1/k_2+1/k_3=1$. Non-trivial hypergeometric solutions of $H_1$
are incomplete elliptic integrals, see Section \ref{ellints}. The degree $d$ is
unbounded, different transformations with the same branching pattern are
possible. Most interesting transformations pull-back the equation $H_1$ into
itself, so that $H_2=H_1$; these transformations come from endomorphisms
of the corresponding elliptic curve.%
\item $1/k_1+1/k_2+1/k_3<1$. Here we have transformations of hyperbolic
hypergeometric functions, see Section \ref{otherat}. The list of these
transformations is finite, the maximal degree of their coverings is 24.
Existence of some of these transformations is shown in \cite{hodgkins1},
\cite{beukers}, \cite{ andkitaev}.
\end{itemize}

The degree of transformations is determined by formula (\ref{trareas}), except in the case 
of incomplete elliptic integrals. If all local exponent differences are real numbers in the
interval $(0,1]$, the covering $\varphi:\PP_x^1\to\PP^1_z$ is defined over
$\RR$ and it branches only above $\{0,1,\infty\}\subset\PP_z^1$, then it
induces a tessellation of the Schwarz triangle for $H_2$ into Schwarz
triangles for $H_1$, as outlined in \cite{hodgkins1,beukers} or
\cite[Section 2]{thyperbolic}. Recall that a {\em Schwarz triangle} for a
hypergeometric equation is the image of the upper half-plane under a Schwarz
map for the equation. The described tessellation is called {\em Coxeter
decomposition}. If it exists, formula (\ref{trareas}) can be interpreted
nicely in terms of areas of the Schwarz triangles for $H_1$ and $H_2$ in the
spherical or hyperbolic metric. Out of the classical transformations, only
the cubic transformation with the branching pattern $3=3=1+1+1$ does not
allow a Coxeter decomposition; see formula (\ref{cubic3}) below.

The following sections form an overview of algebraic transformations for
different types of Gauss hypergeometric functions. We also mention some
three-term identities with Gauss hypergeometric functions. Non-classical
cases are considered %(or will be considered) 
more thoroughly in other articles \cite{tdihedral}, \cite{talggaus}, 
\cite{telliptici}, \cite{thyperbolic}.

\section{Classical transformations}
\label{clalgtr}

Formally, Euler's and Pfaff's
fractional-linear transformations \cite[Theorem 2.2.5]{specfaar}
\begin{eqnarray} \label{flinear1}
\hpg{2}{1}{a,\,b\,}{c}{\,z} & = & (1-z)^{-a}\;
\hpg{2}{1}{a,\,c-b\,}{c}{\frac{z}{z-1}} \\
\label{flinear2} &=& (1-z)^{-b}\;\,
\hpg{2}{1}{c-a,\,b\,}{c}{\frac{z}{z-1}}\\ \label{flinear3}
 & = & (1-z)^{c-a-b}\;\hpg{2}{1}{c-a,\,c-b}{c}{\,z}.
\end{eqnarray}
can be considered as pull-back transformations of degree 1. These are the only
transformations without restrictions on the parameters (or local exponent differences)
of a hypergeometric function under transformation. In a geometrical sense, they 
permute the local exponents at $z=1$ and $z=\infty$.
In general, permutation of the singular points $z=0$, $z=1$, $z=\infty$
and local exponents at them gives 24 Kummer's hypergeometric series solutions
to the same hypergeometric differential equation. Any three hypergeometric solutions
are linearly related, of course.
% Two term relations come from...

To present other classical and non-classical transformations, we
introduce the following notation. Let $(p_1,q_1,r_1)\stackrel{d}{\longleftarrow}(p_2,q_2,r_2)$ schematically denote a pull-back transformation of degree $d$, which transforms a
hypergeometric equation with the local exponent differences $p_1,q_1,r_1$ to
a hypergeometric equation with the local exponent differences $p_2,q_2,r_2$.
The order of local exponents in a triple is irrelevant. 
Note that the arrow follows the direction of the covering \mbox{$\varphi:\PP_x^1\to\PP^1_z$}.
%To get two term relations, one has to place  the point $x=0$ above the 
% or can be adjusted by fractional-linear transformations.

The list of classical transformations comes from the data of Table \ref{clasfig}.
Here is the list of classical transformations with indecomposable $\varphi$,
up to Euler's and Pfaff's fractional-linear transformations and the conversion of
Lemma \ref{basislem}. 
\begin{itemize}
\item $(1/2,\,p,\,q)\stackrel{2}{\longleftarrow}(2p,\,q,\,q)$. These are classical
quadratic transformations. All two-term quadratic transformations of
hypergeometric functions can be obtained by composing (\ref{quadr1}) or
(\ref{quadr2}) with Euler's and Pfaff's transformations. 
%Like we just mentioned, three-term quadratic transformations follow from two-term
%quadratic transformations and three-term fractional-linear transformations.
An example of a three-term relation under a quadratic transformation is the following (see
also Remark \ref{analcont} below, and \cite[2.11(3)]{bateman}):
\begin{eqnarray} \label{quadr3}
\hpg{2}{1}{a,\,b}{\!\frac{a+b+1}{2}}{\,x} & = &
\frac{\Gamma(\frac{1}{2})\,\Gamma(\frac{a+b+1}{2})}
{\Gamma(\frac{a+1}{2})\,\Gamma(\frac{b+1}{2})}\,
\hpg{2}{1}{\frac{a}{2},\,\frac{b}{2}\,}{\frac{1}{2}}{(1\!-\!2x)^2} \nonumber
\\ & & \hspace{-10pt} +(1\!-\!2x)\,\frac{\Gamma(-\frac{1}{2})
\Gamma(\frac{a+b+1}{2})} {\Gamma(\frac{a}{2})\,\Gamma(\frac{b}{2})}\,
\hpg{2}{1}{\!\frac{a+1}{2},\,\frac{b+1}{2}}{\frac{3}{2}}{(1\!-\!2x)^2}.
\end{eqnarray}
\item $(1/2,\,1/3,\,p)\stackrel{3}{\longleftarrow}(1/2,\,p,\,2p)$. These are
well-known Goursat's cubic transformations. Two-term transformations follow
from the following three formulas, along with Euller's and Pfaff's
transformations and application of Lemma \ref{basislem} to (\ref{fcomp6b}):
\begin{eqnarray}
\hpg{2}{1}{a,\,\frac{2a+1}{2}}{\frac{4a+2}{3}}{\,x\,} & \equal &
\left(1\!-\!\frac{3x}{4}\right)^{-a}\,
\hpg{2}{1}{\frac{a}{3},\,\frac{a+1}{3}}{\frac{4a+5}{6}}
{\frac{27\,x^2\,(1\!-\!x)}{(4-3x)^3}},\\
\label{fcomp6a} \hpg{2}{1}{a,\,\frac{2a+1}{2}}{\frac{4a+5}{6}}{\,x\,} &
\equal & \left(1+3x\right)^{-a}\,
\hpg{2}{1}{\frac{a}{3},\,\frac{a+1}{3}}{\frac{4a+5}{6}}
{\frac{27\,x\,(1\!-\!x)^2}{(1+3x)^3}},\\
\label{fcomp6b}\hpg{2}{1}{a,\,\frac{2a+1}{6}}{\frac{1}{2}}{\,x\,} & \equal &
\left(1+\frac{x}{3}\right)^{-a}\,
\hpg{2}{1}{\frac{a}{3},\,\frac{a+1}{3}}{\frac{1}{2}}
{\frac{x\,(9\!-\!x)^2}{(3+x)^3}}.
\end{eqnarray}
\item $(1/3,\,1/3,\,p)\stackrel{3}{\longleftarrow}(p,\,p,\,p)$ These are
less-known cubic transformations. Let $\omega$
denote a primitive cubic root of unity, so $\omega^2+\omega+1=0$. %Here $p=(2a-1)/3$.
Since singular points of the transformed equation are all the same, there is
only one two-term formula (up to changing the parameter):
\begin{equation} \label{cubic3}
\hpg{2}{1}{a,\,\frac{a+1}{3}}{\frac{2a+2}{3}}{\,x}=\left(1+\omega^2
x\right)^{-a}\,\hpg{2}{1}{\frac{a}{3},\,\frac{a+1}{3}}{\frac{2a+2}{3}}
{\frac{3(2\omega\!+\!1)\;x(x-1)}{(x+\omega)^3}}.
\end{equation}
A three-term formula is the following (see also \cite[2.11(38)]{bateman}):
%holds if $|x|<1$, $\mbox{Im}(x)>0$:
% and $\mbox{Re}(a)<\frac{1}{2}$.
\begin{eqnarray}
\hpg{2}{1}{\!a,\frac{a+1}{3}}{\frac{2a+2}{3}}{x} & \equal &
3^{a-1}\,(1\!+\!\omega^2x)^{-a}\left[
\frac{\Gamma(\frac{2a+2}{3})\Gamma(\frac{a}{3})}
{\Gamma(\frac{2}{3})\,\Gamma(a)}\,
\hpg{2}{1}{\!\frac{a}{3},\frac{a+1}{3}}{\frac{2}{3}}
{\left(\frac{x\!+\!\omega^2}{x\!+\!\omega}\right)^3}
\right. \nonumber\\ & & \hspace{2pt} \left. %+(\omega\!+\!1)\,x\,
-\frac{1\!+\!\omega x}{1\!+\!\omega^2x}\,
\frac{\Gamma(\frac{2a+2}{3})\,\Gamma(\frac{a+2}{3})}
{\Gamma(\frac{4}{3})\,\Gamma(a)}\,
\hpg{2}{1}{\!\frac{a+1}{3},\frac{a+2}{3}}{\frac{4}{3}}
{\left(\frac{x\!+\!\omega^2}{x\!+\!\omega}\right)^3} \right].
\end{eqnarray}
\item $(1/2,\,1/3,\,p)\stackrel{4}{\longleftarrow}(1/3,\,p,\,3p)$. These are
indecomposable Goursat's transformations of degree 4. Two-term
transformations follow from the following three formulas, if we compose them
with Euller's and Pfaff's transformations and apply Lemma \ref{basislem} to
(\ref{tr4z}):
\begin{eqnarray}
\hpg{2}{1}{\frac{4a}{3},\,\frac{4a+1}{3}}{\frac{4a+1}{2}}{x} & \equal &
\left(1\!-\!\frac{8x}{9}\right)^{-a}\,
\hpg{2}{1}{\frac{a}{3},\,\frac{a+1}{3}}{\frac{4a+5}{6}}
{\frac{64\,x^3\,(1\!-\!x)}{(9-8x)^3}},\\
\hpg{2}{1}{\frac{4a}{3},\,\frac{4a+1}{3}}{\frac{4a+5}{6}}{x} & \equal &
\left(1+8x\right)^{-a}\,
\hpg{2}{1}{\frac{a}{3},\,\frac{a+1}{3}}{\frac{4a+5}{6}}
{\frac{64\,x\,(1\!-\!x)^3}{(1+8x)^3}},\\
\label{tr4z}\hpg{2}{1}{\frac{4a}{3},\,\frac{4a+1}{6}}{\frac{2}{3}}{x} &
\equal & \left(1-x\right)^{-a}\,
\hpg{2}{1}{\frac{a}{3},\,\frac{3-2a}{6}}{\frac{2}{3}}
{\frac{-x\,(8+x)^3}{64(1-x)^3}}.
\end{eqnarray}
\end{itemize}
As recorded in Table \ref{clasfig}, there are four ways to compose quadratic
and cubic transformations to higher degree transformations of hypergeometric
functions. This gives three different pull-back transformations of degree 4
and 6. The composition transformations can be schematically represented as
follows:
\begin{eqnarray*}
&(1/2,\,1/4,\,p)\stackrel{2}{\longleftarrow}(1/2,\,p,\,p)
\stackrel{2}{\longleftarrow}(p,\,p,\,2p),\\
&(1/2,\,1/3,\,p)\stackrel{3}{\longleftarrow}(1/2,\,p,\,2p)
\stackrel{2}{\longleftarrow}(p,\,p,\,4p),\\ \label{comp6a}
&(1/2,\,1/3,\,p)\stackrel{3}{\longleftarrow}(1/2,\,p,\,2p)
\stackrel{2}{\longleftarrow}(2p,\,2p,\,2p),\\ \label{comp6b}
&(1/2,\,1/3,\,p)\stackrel{2}{\longleftarrow}(1/3,\,1/3,\,2p)
\stackrel{3}{\longleftarrow}(2p,\,2p,\,2p).
\end{eqnarray*}
The last two compositions should produce the same covering, since
computations show that the pull-back
$(1/2,\,1/3,\,p)\stackrel{6}{\longleftarrow}(2p,\,2p,\,2p)$ is unique up to
fractional-linear transformations; see \cite[Section 3]{thyperbolic}.
Indeed, one may check that the identity
\begin{equation} \label{eq:trd6}
\hpg{2}{1}{2a,\,\frac{2a+1}{3}}{\frac{4a+2}{3}}{\,x}=(1-x+x^2)^{-a}
\;\hpg{2}{1}{\!\frac{a}{3},\,\frac{a+1}{3}}{\frac{2a+5}{6}}
{\frac{27}{4}\frac{x^2\,(x\!-\!1)^2}{(x^2\!-\!x\!+\!1)^3}}
\end{equation}
is a composition of (\ref{quadr2}) and (\ref{fcomp6a}), and also a
composition of (\ref{cubic3}), (\ref{quadr2}) and (\ref{flinear1}). Note
that these two compositions use different types of cubic transformations.

\section{Hypergeometric equations with two singularities}
\label{logarithms}

Here we outline transformations of hypergeometric equations with two
relevant singularities; their monodromy group is abelian. The explicit
classification scheme of Section \ref{clscheme} refers to this case three
times. These equations form a special sample of degenerate hypergeometric
equations \cite{degeneratehpg}. For the degenerate cases, not all
usual hypergeometric  formulas for fractional-linear transformations or
other classical algebraic transformations may hold, since the structure of
24 Kummer's solutions degenerates; see \cite[Table 1]{degeneratehpg}.
Here we consider only the new case of pull-back transformations of 
the hypergeometric equations with the cyclic monodromy group.

If a Fuchsian equation has the local exponent difference 1 at some point,
that point can be a non-singular point, an irrelevant singularity or a
logarithmic point. Here is how the logarithmic case is distinguished for
hypergeometric equations.
\begin{lemma} \label{logpoint1}
Consider hypergeometric equation $(\ref{hpgde})$, and let
$P\in\{0,1,\infty\}$. Suppose that the local exponent difference at $S$ is
equal to $1$. Then the point $S$ is logarithmic if and only if (absolute
values of) the two local exponent differences at the other two points of the
set $\{0,1,\infty\}$ are not equal.
\end{lemma}
{\bf Proof.} Because of fractional-linear transformations, we may assume
that $S$ is the point $z=0$, and the local exponents there are 0 and 1.
Therefore $C=0$. Then the point $z=0$ is either a non-singular point or a
logarithmic point. It is non-singular if and only if $A\,B=0$. If $B=0$,
then local exponent differences at $z=1$ and $z=\infty$ are both equal to
$A$.
\qed\\

This lemma implies that a hypergeometric equation has (at most) two relevant
singularities if and only if the local exponent difference at one of the
three points $z=0$, $z=1$, $z=\infty$ is 1, and the local exponent
differences at the other two points are equal. After applying a suitable
fractional-linear transformation to this situation we may assume that the
point $z=0$ is non-singular. Like in the proof of Lemma \ref{logpoint1}, we
have $C=0$ and we may take $B=0$. Then we are either in the case $n=m=0$ of
\cite[Section 7 or 8]{degeneratehpg}, or in the case $n=m=\ell=0$ of
\cite[Section 9]{degeneratehpg}. Most of the 24 Kummer's solutions have to
be interpreted either as the constant 1 or the power function $(1-z)^{-a}$.
The only interesting hypergeometric function (up to Euler's and Pfaff's
transformations) is the following:
\begin{equation} \label{reghpg}
\hpg{2}{1}{1+a,\,1}{2}{\,z\,}=\left\{ \begin{array}{cl} \displaystyle
\frac{(1-z)^{-a}-1}{a\,z}, & \mbox{if } a\neq 0, \vspace{3pt}\\
\displaystyle -\frac1z\,\log(1-z), & \mbox{if } a=0. \end{array} \right.
\end{equation}

For general $a$, pull-back transformation (\ref{algtransf}) of the
considered hypergeometric equation to a hypergeometric equation branches
only above the points $z=1$ and $z=\infty$. Indeed, if the covering
$\varphi:\PP_x^1\to\PP^1_z$ branches above other point, then these
branching points would be singular by part 2 of Lemma \ref{genrami}, and
there would be at least 3 singular points above $\{1,\infty\}\subset\PP^1_z$ %(essentially 
by part 1 of Lemma \ref{genrami2}. To keep the number of
singular points down to 3, the covering $\varphi$ should branch only
above $\{1,\infty\}$. Up to fractional-linear transformations on $\PP^1_x$,
these coverings have the form $(1-z) \longmapsto (1-x)^d$, or
\begin{equation} \label{cyclictr}
z\longmapsto x\,\phi_{d-1}(x), \qquad \mbox{where} \quad
\phi_{d-1}(x)=\frac{1-(1-x)^d}{x}.
\end{equation}
Note that $\phi_{d-1}(x)$ is a polynomial of degree $d-1$. A corresponding
hypergeometric identity is
\begin{equation}
\hpg{2}{1}{1+da,\,1}{2}{\,x\,}=\frac{\phi_{d-1}(x)}{d}\;
\hpg{2}{1}{1+a,\,1}{2}{x\,\phi_{d-1}(x)}.
\end{equation}
This transformation is obvious from the explicit expressions in
(\ref{reghpg}). 

Formally, we additionally have a continuous family $z\mapsto 1-\beta+\beta z$ of
fractional-linear pull-back transformations which fix the two points $z=1$ and $z=\infty$.
However, they do not give interesting hypergeometric
identities since Kummer's series at those two points are trivial.

If $|a|=1/k$ for an integer $k>1$, there are more pull-back transformations 
of hypergeometric equations with the local exponent differences $(1,a,a)$.  
In this case, the monodromy group is a finite cyclic group, of order $k$.
Pull-backed equations will have a cyclic monodromy group as well, possibly
of smaller order. On the other hand, the mentioned Klein's theorem \cite{klein77}
%(on second order differential equations with a finite monodromy group)
implies that any hypergeometric equation with a cyclic monodromy group of order $k$
is a pull-back of a hypergeometric equation with the local exponent differences
$(1,1/k,1/k)$. These pull-back transformations can be easily computed from
explicit terminating solutions of the target differential equation.
According to \cite[Section 7]{degeneratehpg}, a general hypergeometric equation
with  a completely reducible (but non-trivial) monodromy representation has the local exponents
$(m+n+1,a,a+n-m)$, where $a\not\in\ZZ$ and $n,m\in\ZZ$ are non-negative.
A basis of terminating solutions is 
\begin{equation} \label{eq:credb}
\hpg21{-n,a-m}{-m-n}{\,z\,},\qquad (1-z)^{-a}\,\hpg21{-m,-a-n}{-m-n}{\,z\,}.
\end{equation}
The monodromy group is finite cyclic if $a=\ell/k$ with co-prime positive $k,\ell\in\ZZ$.
The terminating solutions can be written as terminating hypergeometric series at $z=1$ as well:
\begin{eqnarray*}
\hpg21{-n,a-m}{-m-n}{\,z\,}= \frac{(1+a)_n\,m!}{(m+n)!}\,
\hpg21{-n,a-m}{1+a}{1-z}, \qquad\mbox{etc.}
\end{eqnarray*}
The quotient of two solutions in (\ref{eq:credb}) defines a {\em Schwarz map} for the 
hypergeometric equation. In the simplest case $n=m=0$, $a=1/k$, the Schwartz map
is just $(1-z)^{1/k}$. Klein's pull-back transformation for 
$(1,1/k,1/k)\stackrel{d}{\longleftarrow}(m+n+1,\ell/k,\ell/k+n-m)$
is obtained from identification of the two Schwarz maps.
The pull-back covering is defined by
 \begin{equation} \label{eq:padetr}
(1-z)\longmapsto (1-x)^{\ell}\,
\hpg21{-n,\,\ell/k-m}{-m-n}{x}^k\left/
\hpg21{-m,-\ell/k-n}{-m-n}{x}^k\right..
\end{equation}
The Schwarz maps (or pairs of hypergeometric solutions) 
are identified here by the corresponding local exponents at $x=1$ (placed above $z=1$)
and the same value at $x=0$ (placed above $x=0$).
The degree of the transformation is equal to $\max(nk+\ell,mk)$, by formula (\ref{trareas}) as well.
Besides, $z\longmapsto O(x^{n+m+1})$ at $x=0$ by the required branching pattern. 
In particular, 
\begin{eqnarray}
\hpg21{-m,-\ell/k-n}{-m-n}{x}\left/
\hpg21{-n,\,\ell/k-m}{-m-n}{x}\right.=(1-x)^{\ell/k}+O(x^{n+m+1})
%\qquad\mbox{at}\quad x=0
\end{eqnarray}
at $x=0$, hence the quotient of two hypergeometric polynomials is
the {\em Pad\'e approximation} of $(1-x)^{\ell/k}$ of precise degree $(m,n)$.
For example, the Pade approximation of $\sqrt{1-x}$ of degree $(1,1)$ is
$(4-x)/(4-3x)$. Hence the following pullback must give a transformation
$(1,1/2,1/2)\stackrel{3}{\longleftarrow}(3,1/2,1/2)$:
\[
1-z\longmapsto \frac{(1-x)\,(x-4)^2}{(3x-4)^2}.
\]
A corresponding hypergeometric identity is
\begin{equation}
\hpg{2}{1}{3/2,\,2\,}{4}{\,x\,}=\frac{4}{4-3x}\;
\hpg{2}{1}{1/2,\,1\,}{2}{\frac{x^3}{(3x-4)^2}}.
\end{equation}

Transformation (\ref{eq:padetr}) is Klein's pull-back transformation if $\gcd(k,\ell)=1$.
Otherwise the transformed hypergeometric equation has a smaller monodromy group.
These transformations must factor via (\ref{cyclictr}) with $d=\gcd(k,\ell)$,
and Klein's transformation between equations with the smaller monodromy group.
Even $\ell/k\in\ZZ$ can be allowed if the transformed equation has no logarithmic points.
The condition for that is $\ell/k>m$; see \cite[Corollary 2.3 part (2)]{degeneratehpg}.
Under this condition, one may even allow $k=1$ and consider transformations
$(1,1,1)\stackrel{\;\ell+n}{\longleftarrow}(m+n+1,\ell,\ell+n-m)$. All hypergeometric equations
with the trivial monodromy group can be obtained in this way, by Klein's theorem.
Solutions of these hypergeometric equations 
are analyzed  in \cite[Section 8]{degeneratehpg}. A hypergeometric equation with the local exponent differences $(1,1,1)$ can be transformed to $y''=0$ by fractional-linear transformations.
We underscore that transformation (\ref{eq:padetr}) specializes nicely even for $k=1$
if only logarithmic solutions are not involved; the corresponding two-term hypergeometric
identities are trivial.

\vspace{2pt}
\begin{remark} \label{analcont} \rm
Algebraic transformations of Gauss hypergeometric functions often hold only
in some part of the complex plane, even after standard analytic continuation. For
example, formula (\ref{quadr1}) is obviously false at $x=1$. Formula (\ref{quadr1}) 
holds when $\mbox{Re}(x)<1/2$, as the standard $z$-cut $(1,\infty)$ is mapped into the
line $\mbox{Re}(x)=1/2$ under the transformation $z=4x(1-x)$.

An extreme example of this kind is the following transformation of a hypergeometric function to
a rational function:
\begin{equation} \label{ratff}
\hpg{2}{1}{1/2,\,1\,}{2}{-\frac{4\,x^3\,(x-1)^2\,(x+2)}
{(3x-2)^2}}=\frac{2-3x}{(1-x)^2\,(x+2)}.
\end{equation}
This identity holds in a neighborhood of $x=0$, but it certainly does not
hold around $x=1$ or $x=-2$. Apparently, standard cuts for analytic
continuation for the hypergeometric function isolate the three points $x=0$,
$x=1$, $x=-2$. Note that $\hpg{2}{1}{\!1/2,\,1}{2}{z}=(2-2\sqrt{1-z})/z$ is a two-valued
algebraic function on $\PP_z^1$. % by quadratic transformation (\ref{quadr1}).
Its composition in (\ref{ratff}) with the degree 6 rational function 
apparently consists of two disjoint branches. The second branch is the
rational function $(3x-2)/x^3$, which is the correct evaluation of the
left-hand side of (\ref{ratff}) around the points $x=1$, $x=-2$ (check the
power series.)

Many identities like (\ref{ratff}) can be produced for hypergeometric
functions of this section with $1/a\in\ZZ$. The pull-backed hypergeometric
equations should be Fuchsian equations with the trivial monodromy group.
More generally, any algebraic hypergeometric function can be pull-backed to
a rational function. Other algebraic hypergeometric functions are considered
in the following two sections.

Three-term hypergeometric identities may also have limited region of
validity. But it may happen that branch cuts of two hypergeometric terms
cancel each other in a three-term identity. For example, standard branch
cuts for the hypergeometric functions on the right-hand side of
(\ref{quadr3}) are the intervals $[1,\infty)$ and $(-\infty,0]$ on the real
line. But identity (\ref{quadr3}) is valid on $\CC\setminus[1,\infty)$, if
we agree to evaluate the right-hand consistently on the interval
$(-\infty,0]$: either using analytic continuation of both terms from the
upper half-plane, or from the lower half-plane.
\end{remark}

\section{Dihedral functions} \label{dihedrals}

Hypergeometric equations with (infinite or finite) dihedral monodromy group
are characterized by the property that two local exponent differences are
rational numbers with the denominator 2. By a quadratic pull-back
transformation, these equations can be transformed to Fuchsian equations
with at most 4 singularities and with a cyclic monodromy group. Explicit
expressions and transformations for these functions are considered thoroughly
in \cite{tdihedral}, \cite{dihedraltr}. Here we look at transformations 
of hypergeometric equations  which have two local exponent differences equal to $1/2$. 
The explicit classification scheme of Section \ref{clscheme} refers to this case twice.

The starting hypergeometric equation for new transformations has
the local exponent differences $(1/2,1/2,a)$. Hypergeometric solutions of such an equation
can be written explicitly. In particular, quadratic transformation (\ref{quadr1}) with $b=a+1$
implies
%\begin{equation} \label{fpdihedr}
%\hpg{2}{1}{\frac{a}{2},\,\frac{a+1}{2}}{a+1}{\,4x(1-x)\,}=
%\hpg{2}{1}{a,\,a+1}{a+1}{\,x}=(1-x)^{-a}.
%\end{equation}
\begin{eqnarray} \label{dihedr1}
\hpg{2}{1}{\frac{a}{2},\,\frac{a+1}{2}\,}{a+1}{\,z\,} & = &
\left(\frac{1+\sqrt{1-z}}{2} \right)^{-a}.
\end{eqnarray}
Other explicit formulas are
\begin{eqnarray} \label{dihedr2} 
%\hpg{2}{1}{\frac{a}{2},\,\frac{a+1}{2}\,}{a+1}{\,z\,} & = &
%\left(\frac{1+\sqrt{1-z}}{2} \right)^{-a},\\ \label{dihedr2}
\hpg{2}{1}{\frac{a}{2},\,\frac{a+1}{2}\,}{\frac{1}{2}}{\,z\,} & = &
\frac{(1-\sqrt{z})^{-a}+(1+\sqrt{z})^{-a}}{2},\\ \label{dihedr3}
\hpg{2}{1}{\!\frac{a+1}{2},\,\frac{a+2}{2}}{\frac{3}{2}}{\,z\,} & = &
\left\{\begin{array}{cl} \displaystyle
\frac{(1-\sqrt{z})^{-a}-(1+\sqrt{z})^{-a}}{2\,a\,\sqrt{z}}, &\mbox{if}\
a\neq 0,\\ \displaystyle
\frac1{2\sqrt{z}}\,\log\frac{1+\sqrt{z}}{1-\sqrt{z}}, & \mbox{if}\ a=0.
\end{array}\right.
\end{eqnarray}
General dihedral Gauss hypergeometric functions are contiguous to these $\hpgo21$ functions.
As shown in \cite{tdihedral},
explicit expressions for them can be given in terms of terminating Appell's $F_2$ or $F_3$ series.
For example, generalizations of (\ref{dihedr1})--(\ref{dihedr2}) are
\begin{eqnarray}  \label{eq:diha} 
\hpg{2}{1}{\frac{a}{2},\,\frac{a+1}{2}+\ell}{a+k+\ell+1}{\,1-z\,} &\equal&
z^{k/2} \left(\frac{1+\sqrt{z}}{2} \right)^{-a-k-\ell}\times\nonumber\\
&&\app3{k+1,\ell+1; -k,-\ell}{a+k+\ell+1}{\frac{\sqrt{z}-1}{2\sqrt{z}},\frac{1-\sqrt{z}}2},\\
\label{eq:dih12} %\hspace{-9pt}
\frac{\left(\frac{a+1}2\right)_n}{\left(\frac12\right)_n}\,
\hpg{2}{1}{\frac{a}{2},\frac{a+1}{2}+n}{\frac{1}{2}-m}{\,z}\! &\equal& \frac{(1+\sqrt{z})^{-a}}2
\app2{a; -m,-n}{-2m,-2n}{\frac{2\sqrt{z}}{1+\sqrt{z}},\frac2{1+\sqrt{z}}}\nonumber\\
&&+\frac{(1-\sqrt{z})^{-a}}2
\app2{a; -m,-n}{-2m,-2n}{\frac{2\sqrt{z}}{\sqrt{z}-1},\frac2{1-\sqrt{z}}}.
\end{eqnarray}
Here $m,n$ are assumed to be non-negative integers.

For general $a$, there are two types of transformations:
\begin{itemize}
\item $(1/2,1/2,a)\stackrel{d}{\longleftarrow}(1/2,\,1/2,\,da)$. These are
the only transformations to a dihedral monodromy group as well, as there is a singularity 
above the point with the local exponent difference $a$. 
Identification of explicit Schwarz maps gives the following recipe for
computing the pull-back coverings $\varphi:\PP_x^1\to\PP_z^1$. 
Expand $(1+\sqrt{x})^d$ in the form $\theta_1(x)+\theta_2(x)\sqrt{x}$ with
$\theta_1(x),\theta_2(x)\in\CC[x]$. Then $\varphi(x)=x\,\theta^2_2(x)/\theta_1^2(x)$
gives a pull-back transformation of dihedral hypergeometric equations. Explicitly,
\begin{eqnarray*} \label{thetasd}
\theta_1(x)=&\displaystyle\sum_{k=0}^{\lfloor d/2 \rfloor}\,{d\choose
2k}\;x^k & = \hpg21{-\frac{d}2,-\frac{d-1}2}{\frac12}{\;x\,}, \\
\theta_2(x)=&\displaystyle\!\!\sum_{k=0}^{\lfloor (d-1)/2 \rfloor}\!
 {d\choose 2k\!+\!1}\,x^k&=d\;\hpg21{\!-\frac{d-1}2,-\frac{d-2}2}{\frac32}{\,x}.
\end{eqnarray*}
A particular transformation of hypergeometric functions is the following:
\begin{equation} \label{eq:gdihtr}
\hpg{2}{1}{\frac{da}{2},\,\frac{da+1}{2}}{\frac{1}{2}}{\,x}=\theta_1(x)^{-a}
\,\hpg{2}{1}{\frac{a}{2},\,\frac{a+1}{2}}{\frac{1}{2}}
{\frac{x\,\theta_2(x)^2}{\theta_1(x)^2}}.
\end{equation}
It is instructive to check this transformation using (\ref{dihedr2}).
Other transformations from the same pull-back covering are given in \cite[Section 6]{tdihedral}.
Particularly interesting are the following formulas; they hold for odd or even $d$, respectively:
\begin{eqnarray*}
\hpg{2}{1}{\frac{da}{2},-\frac{da}{2}}{\frac12}{\,x} &\equal&
\hpg{2}{1}{\frac{a}{2},-\frac{a}{2}}{\frac12}{d^2x\,\hpg21{\frac{1-d}2,\frac{1+d}2}{3/2}{x}^{\!2\,}},\\
\hpg{2}{1}{\frac{da}{2},-\frac{da}{2}}{\frac12}{\,x} &\equal&
\hpg{2}{1}{\frac{a}{2},-\frac{a}{2}}{\frac12}{d^2x(1-x)\hpg21{1-\frac{d}2,1+\frac{d}2}{3/2}{x}^{\!2\,}}.
\end{eqnarray*}
The branching pattern of $\varphi(x)$ is
\begin{eqnarray*}
1+2+2+\ldots+2=d=1+2+2+\ldots+2, && \mbox{if $d$ is odd},\\
1+1+2+2+\ldots+2=d=2+2+\ldots+2,\hspace{17pt} && \mbox{if $d$ is even}.
\end{eqnarray*}
\item $(1/2,1/2,a)\stackrel{2\ell}{\longleftarrow}(1,\ell a,\,\ell\,a)$, and
$d=2\ell$ is even. These are transformations to hypergeometric equations of
Section \ref{logarithms}. They are compositions of the mentioned quadratic
transformation and the transformations 
$(1/2,1/2,a)\stackrel{d}{\longleftarrow}(1/2,1/2,da)$ or 
$(1,a,a)\stackrel{d}{\longleftarrow}(1,d a,d a)$ described above.
\end{itemize}
If $a=1/k$ with $k$ a positive integer, the monodromy group is the
finite dihedral group with $2k$ elements, and hypergeometric solutions are
algebraic. Klein's theorem \cite{klein77} implies that any hypergeometric equation
with a finite dihedral monodromy group is a pull-back from a hypergeometric equation
with the local exponent differences $(1/2,1/2,1/k)$ and the same monodromy group. 
The pull-back transformation can be computed by the similar method:
identification of explicit Schwarz maps, using the mentioned explicit evaluations
with terminating Appell's $F_2$ or $F_3$ series. That leads to expressing a polynomial
in $\sqrt{x}$ in the form \mbox{$\theta_1(x)+\sqrt{x}\theta_2(x)$} as above.
\begin{theorem} \label{th:gfdih}
Let $k,\ell,m,n$ be positive integers, and suppose that $k\ge 2$, $\gcd(k,\ell)=1$. 
Let us denote
\begin{eqnarray*} \label{eq:gdef1}
G(x)= %\left(1+\sqrt{x}\right)^{n}
x^{m/2}\,\app3{m+1,n+1; -m,-n}{1+\ell/k}{\frac{\sqrt{x}+1}{2\sqrt{x}},\frac{1+\sqrt{x}}2}.
\end{eqnarray*}
This is a polynomial in $\sqrt{x}$. We can write 
\begin{eqnarray*} \label{eq:gmpower}
\left(1+\sqrt{x}\right)^{\ell}G(x)^k=\Theta_1(x)+x^{m+\frac12}\,\Theta_2(x),
\end{eqnarray*}
so that $\Theta_1(x)$ and $\Theta_2(x)$ are polynomials in $x$. 
Then the rational function $\Phi(x)=x^{2m+1}\,\Theta_2(x)^2/\Theta_1(x)^2$
defines Klein's pull-back covering 
\mbox{$(1/2,1/2,1/k)\stackrel{d}{\longleftarrow}(m+1/2,n+1/2,\ell/k)$}.
The degree $d$ of this rational function %$\Phi(x)$ 
is equal to $(m+n)k+\ell$.
\end{theorem}
\proof This is Theorem 5.1 in \cite{dihedraltr}. %Theorem 7.1 in \cite{tdihedral}.
\qed\\

\noindent
The condition $\gcd(k,\ell)=1$ can be replaced by the weaker condition $\ell/k\not\in\ZZ$,
but then the transformed hypergeometric equation has a smaller dihedral
monodromy group, and it factors via the transformation in (\ref{eq:gdihtr}) with $d=\gcd(k,\ell)$.
Even more, $\ell/k\in\ZZ$ %or even $k=1$ 
can be allowed, if the transformed equation has no logarithmic solutions. 
Sufficient and necessary conditions for that are given in \cite[Theorem 2.1]{tdihedral}. 
The branching pattern for all these coverings has the following pattern: 
\begin{itemize}
\item Above the two points with the local exponent difference $1/2$,
there are two points with the branching orders  $2m+1$, $2n+1$,
and the remaining points are simple branching points.
\item Above the point with the local exponent difference $1/k$,
there is one point with the ramification order $\ell$, 
and $m+n$ points with the ramification order $k$.
\end{itemize}
Any covering $(1/2,1/2,1/k)\stackrel{d}{\longleftarrow}(m+1/2,n+1/2,\ell/k)$ is unique
up to fractional-linear transformations, as Schwarz maps are identified 
uniquely\footnote{If only $k$ does not divide $\ell$. See \cite[Subsection 5.3]{dihedraltr} for a counterexample with $(m,n,k,\ell)=(0,0,5,10)$.}. 
Transformations from the local exponent differences $(1/2,1/2,1/k)$ 
to hypergeometric equations with finite cyclic monodromy groups are 
either the mentioned degeneration $\ell/k\in\ZZ$, or compositions
with the quadratic transformation $(1/2,1/2,1/k)\stackrel{2}{\longleftarrow}(1,1/k,1/k)$. 
Other transformations involving dihedral Gauss hypergeometric functions
are special cases of classical transformations.

For the purposes of Theorem \ref{th:gfdih}, the function $G(x)$ can be alternatively
defined as follows:
\begin{equation}
\left(1+\sqrt{x}\right)^{(m+n)+\ell}\,\app2{-\ell/k-m-n;-m,-n}{-2m,-2n}
{\frac{2\sqrt{x}}{1+\sqrt{x}},\frac2{1+\sqrt{x}}}^k.
\end{equation}
The two definitions differ by a constant multiple. The $F_2$ and $F_3$ sums are
related by reversing the order of summation in both directions in the rectangular sums,
as noted in \cite{tdihedral}.

For an example, consider the case $n=1,m=0,\ell=1$  of Theorem \ref{th:gfdih}.  
To compute the transformation
$(1/2,1/2,1/k)\stackrel{\;k+1}{\longleftarrow}(1/2,\,3/2,\,1/k)$ we need to expand
\[
\left(1+\sqrt{x}\right)\left(1-\frac{\sqrt{x}}{k}\right)^k=
\theta_3\!\left(x\right)+x^{3/2}\,\theta_4\!\left(x\right).
\]
Straightforward computation shows that
\begin{eqnarray} \label{theta34}
\theta_3(x)=\hpg21{-\frac{k}2,-\frac{k+1}2}{-1/2}{\frac{x}{k^2}},\qquad
\theta_4(x)=\frac{k^2-1}{3k^2}\,\hpg21{-\frac{k-2}2,-\frac{k-3}2}{5/2}{\frac{x}{k^2}}.
\end{eqnarray}
A transformation of hypergeometric functions is
\begin{equation}
\hpg21{-\frac{1}{2k},-\frac12-\frac{1}{2k}}{-\frac12}{\,x}=
\theta_3(x)^{1/k}
\hpg21{-\frac{1}{2k},\frac12-\frac{1}{2k}}{\frac12}{\frac{x^3\,\theta_4(x)^2}{\theta_3(x)^2}}.
\end{equation}
On the other hand,
\begin{eqnarray*}
\hpg21{-\frac{1}{2k},\,-\frac12-\frac{1}{2k}}{1-\frac{1}k}{1-z}=\frac{k-\sqrt{z}}{k-1}
\left(\frac{1+\sqrt{z}}2\right)^{1/k}
\end{eqnarray*}
by formula (\ref{eq:diha}). Note that the construction in (\ref{theta34}) breaks down if $k=1$; 
a hypergeometric equation with the local exponent differences $(1/2,\,3/2,\,1)$ has logarithmic solutions.

As computed in \cite[Section 7]{tdihedral}, the polynomials  $\Theta_1(x)$, $\Theta_2(x)$
of Theorem \ref{th:gfdih} in the case $n=1$, $m=0$, $\ell=2$ can be expressed as terminating 
$\hpgo32$ series.

%For general dihedral equations, the rational function $\varphi(x)$ cannot be
%expressed as a quotient of terminating hypergeometric series.
%%This is connected to the fact that pull-back transformation of a dihedral
%%hypergeometric equation to a Fuchsian equation is not hypergeometric, but has 4 singularities.

\section{Algebraic Gauss hypergeometric functions} \label{alggtranfs}

Algebraic Gauss hypergeometric functions form a classical subject of
mathematics. These functions were classified by Schwarz \cite{Schwarz72}.
Recall that a Fuchsian equation has a basis of algebraic solutions if and
only if its monodromy group is finite. Finite {\em projective monodromy groups} for second
order equations are either cyclic, or dihedral, or the tetrahedral group
isomorphic to $A_4$, or the octahedral group isomorphic to $S_4$, or the
icosahedral group isomorphic to $A_5$. An important characterization of
second order Fuchsian equations with finite monodromy group was given by
Klein \cite{klein77,klein78}: all these equations are pull-backs of a few {\em
standard hypergeometric equations} with algebraic solutions. In particular,
this holds for hypergeometric equations with finite monodromy groups. The
corresponding standard equation depends on the projective monodromy group:
\begin{itemize}
\item Second order equations with a cyclic monodromy group are pull-backs of
a hypergeometric equation with the local exponent differences $(1,1/k,1/k)$,
where $k$ is a positive integer. Klein's transformations to general
hypergeometric equations with a cyclic monodromy group are
considered in Section \ref{logarithms} above.%
\item Second order equations with a finite dihedral monodromy group are
pull-backs of a hypergeometric equation with the local exponent differences
$(1/2,1/2,1/k)$, where $k\ge 2$. Klein's transformations to general
hypergeometric equations with a dihedral monodromy group are
considered in Section \ref{dihedrals} above.%
\item Second order equations with the tetrahedral projective monodromy group are
pull-backs of a hypergeometric equation with the local exponent differences
$(1/2,1/3,1/3)$. Hypergeometric equations with this monodromy group are
contiguous to hypergeometric equations with the local exponent differences
$(1/2,1/3,1/3)$ or $(1/3,1/3,2/3)$.  %
\item Second order equations with the octahedral projective monodromy group are
pull-backs of a hypergeometric equation with the local exponent differences
$(1/2,1/3,1/4)$. Hypergeometric equations with this monodromy group are
contiguous to hypergeometric equations with the local exponent differences
$(1/2,1/3,1/4)$ or $(2/3,1/4,1/4)$.%
\item Second order equations with the icosahedral projective monodromy group are
pull-backs of a hypergeometric equation with the local exponent differences
$(1/2,1/3,1/5)$. Hypergeometric equations with this monodromy group are
contiguous to hypergeometric equations with the local exponent differences
$(1/2,1/3,1/5)$, $(1/2,1/3,2/5)$, $(1/2,1/5,2/5)$, $(1/3,1/3,2/5)$,
$(1/3,2/3,1/5)$, $(2/3,1/5,1/5)$, $(1/3,2/5,3/5)$, $(1/3,1/5,3/5)$,
$(1/5,1/5,4/5)$ or $(2/5,2/5,2/5)$.
\end{itemize}
A general algorithm for computation of Klein's coverings is given in \cite{kleinvhw}.
The algorithm is based on finding semi-invariants of the monodromy group by
solving appropriate symmetric powers of the given second order differential equation.
A more effective algorithm specifically for hypergeometric equations with finite monodromy
groups is given in \cite{talggaus}. This algorithm is based on identification of explicit Schwarz
maps for the given %hypergeometric equation 
and the corresponding standard hypergeometric equation.

Klein's pull-back transformations are most interesting among the transformations of algebraic $\hpgo21$
functions. As we will show soon, all other transformations between hypergeometric equations with a finite monodromy group are special cases of classical transformations, expect the series of transformations of Sections \ref{logarithms}, \ref{dihedrals}, and one degree 5 transformation between standard icosahedral and octahedral hypergeometric equations.

First we sketch the algorithm in  \cite{talggaus} for computing Klein's pull-back transformations of hypergeometric equations with the tetrahedral, octahedral or icosahedral 
projective monodromy groups. The mentioned contiguous orbits of hypergeometric functions
determine {\em Schwarz types} of algebraic $\hpgo21$  functions. There is one dihedral,
2 tetrahedral, 2 octahedral and 10 icosahedral Schwarz types.

The algorithm in \cite{talggaus} uses explicit evaluation of algebraic Gauss hypergeometric functions,
called {\em Darboux evaluations}. The geometric idea behind them is to pull-back a hypergeometric equation with a finite monodromy group to a Fuchsian differential equation with a cyclic
monodromy group. Pull-backed hypergeometric solutions can be expressed in terms of radical functions, like in formulas (\ref{dihedr1})--(\ref{eq:dih12}) for dihedral functions.
The minimal degree for these {\em Darboux pull-backs} to a cyclic monodromy group
is 3, 4 or 5 for, respectively, tetrahedral, octahedral and icosahedral differential equations. 
The quadratic transformation $(1/2,1/2,a)\stackrel{2}{\longleftarrow}(1,a,a)$ in Section \ref{dihedrals}
is actually a Darboux pull-back in the dihedral case. Here are a few examples of Darboux evaluations
for larger finite monodromy groups:
\begin{eqnarray*} \label{fptetra1}
\hpg{2}{1}{1/4,-1/12\,}{2/3}{\,\frac{x\,(x+4)^3}{4(2x-1)^3}} &\equal&
\frac1{\left(1-2x\right)^{1/4}},\\
\hpg{2}{1}{1/2,\,-1/6\,}{2/3}{\,\frac{x\,(x+2)^3}{(2x+1)^3}}
&\equal& \frac1{\sqrt{1+2x}},\\
\hpg{2}{1}{7/24,\,-1/24}{3/4}{\frac{108\,x\,(x-1)^4}{(x^2+14x+1)^3}}
&\equal& \frac1{(1+14x+x^2)^{1/8}},\\
\hpg{2}{1}{1/6,-1/6\,}{1/4}{\frac{27\,x\,(x+1)^4}{2(x^2+4x+1)^3}} &\equal&
\frac{\left(1+2x\right)^{1/4}}{\sqrt{1+4x+x^2}},\\
\hpg{2}{1}{\!13/60,-7/60}{3/5}{\frac{1728\,x\,(x^2-11x-1)^5}{(x^4\!+\!228x^3\!+\!494x^2\!-\!228x\!+\!1)^3}}
& \equal &\frac{1-7x}{\big(1\!-\!228x\!+\!494x^2\!+\!228x^3\!+\!x^4\big)^{7/20}},\\
\hpg{2}{1}{7/20,\,-1/20}{4/5}{\frac{64\,x\,(x^2-x-1)^5}{(x^2\!-\!1)\,(x^2+4x-1)^5}}&\equal&
\frac{\left(1+x\right)^{7/20}}{\left(1-x\right)^{1/20}\left(1-4x-x^2\right)^{1/4}}.\hspace{56pt}
\end{eqnarray*}
Some of minimal Darboux pull-back coverings for icosahedral functions are defined
not on $\PP^1_x$, but on a genus 1 curve. For example,
\begin{eqnarray*}
\hpg{2}{1}{8/15,-1/15}{4/5}{\frac{54\;(\xi_1+5x)^3(1-2\xi_1+6x)^5}
{(16x^2\!-\!1)(\xi_1\!-5x)^2(1-2\xi_1\!-14x)^5}} & \equal &
\frac{(1+4x)^{8/15}\,(\xi_1+5x)^{1/6}\,x^{1/15}}
{(1-2\xi_1-14x)^{1/3}\;(\xi_1-3x)^{3/10}},\\
\hpg{2}{1}{7/10,-1/10\,}{4/5}{\,\frac{16\,\xi_2\,(1+x-x^2)^2\,(1-\xi_2)^2}
{(1+\xi_2+2x)(1+\xi_2-2x)^5}} & \equal & \frac{
(1-\xi_2+2x)^{1/15}\,(1-\xi_2)^{3/5}}{(1+\xi_2+2x)^{7/30}\,\sqrt{1+\xi_2-2x}},
\end{eqnarray*}
where
\[
\xi_1=\sqrt{x\;(1+x)\,(1+16x)}\qquad\mbox{and}\qquad\xi_2=\sqrt{x\;(1+x-x^2)}.
\]
These formulas can be checked with a computer algebra package by expanding
both sides in power series in $x$ or $\sqrt{x}$. In \cite{talggaus}, a few of these evaluations
are computed for each Schwarz type of algebraic Gauss hypergeometric functions. 
Using contiguous relations, one can find a Darboux evaluation for any algebraic $\hpgo21$ function.
For comparison, in Section \ref{dihedrals} we used general formulas (with terminating
Appell's $F_2$ or $F_3$ sums) for dihedral $\hpgo21$ functions,
instead of applying contiguous relations.

A suitable ramification pattern for Klein's pull-back covering, say for a transformation $(1/2,1/3,1/4)$
$\stackrel{d}{\longleftarrow}(n+1/2,\,m+1/3,\,\ell+1/4)$ of local exponent differences,
is easy to set up. In the setting of Section \ref{preliminar}, it is convenient to assume that
$x=0$ lies above $z=0$ and assign local exponent differences with the largest denominator 
(say, 4) to these points. Hypergeometric solutions of the given and its standard hypergeometric equations at these points can be {\em a priori} identified 
(up to a constant multiple, at worst) by their local exponents. This gives identification of 
Schwarz maps (for both hypergeometric equations) up to a constant multiple.
The constant multiple can be determined by a separate routine for each Schwarz type.
Elimination of the variables in Darboux evaluations gives an algebraic relation
between the arguments of the given and its standard hypergeometric equations,
which must give Klein's pull-back covering. The degree of Klein's pull-back covering
for a transformation $(1/2,1/3,1/k)\stackrel{d}{\longleftarrow}(e_0,\,e_1,\,e_\infty)$ 
can be computed from (\ref{trareas}):
\begin{equation}
d %=\frac{e_0+e_1+e_\infty-1}{\frac12+\frac13+\frac1m-1}
= \frac{6k}{6-k}(e_0+e_1+e_\infty-1).
\end{equation}

Here is a list of Klein's pull-back coverings computed in \cite{talggaus}.
The triples of local exponent differences for the transformed equation are given on the left.
The standard hypergeometric equations have the local exponent differences
$(1/2,1/3,1/3)$ or $(1/2,1/3,1/5)$. 
\begin{eqnarray*}
(1/2,2/3,2/3): && z=-\frac{x^2(4x-5)^3}{(5x-4)^3},\\
(3/2,1/3,1/3): && z = -\frac{x(x^2-42x-7)^3}{(7x^2+42x-1)^3},\\
(1/2,1/3,4/3): && z = -\frac{x(256x^2-448x+189)^3}{27(28x-27)^3},\\
(1/2,2/3,4/3): && z = \frac{19683x^2(4x-1)^3}{(256x^3-192x^2+21x-4)^3},\\
(1/2,1/3,5/3): && z = -\frac{19683x(128x-125)^3}{(16384x^3-30720x^2+14880x-625)^3},\\
(3/2,1/3,2/3): && z = -\frac{729x(5x^2+14x+125)^3}{(4x^3+15x^2-690x-625)^3},\\
(1/3,2/3,5/3): && z = \frac{4x(256x^3-640x^2+520x-135)^3}{27(x-1)^2(32x-27)^3},\qquad \\
(2/3,2/3,4/3): && z = -\frac{x^2(x-1)^2(16x^2-16x+5)^3}{4(5x^2-5x+1)^3},\\
(2/3,4/3,4/3): && z=-\frac{108\;x^4\;(x-1)^4\;(27x^2-27x+7)^3}{(189x^4-378x^3+301x^2-112x+16)^3},\\
(1/2,2/3,1/5): && z=\frac{x\,(102400x^2-11264x-11)^5}{(180224000x^3+4325376x^2-21252x+1)^3},\\
(1/5,1/5,6/5): && \!z=\frac{108\,x\,(1-x)\,(512x^2-512x+3)^5}
{(1048576x^6\!-\!3145728x^5\!+\!3244032x^4\!-\!1245184x^3\!+\!94848x^2\!+\!3456x\!+\!1)^3}.
\end{eqnarray*}
Once a pull-back covering is known,  hypergeometric identities are easy to derive.
For example, 
\begin{eqnarray} \label{eq:hpg5}
\hpg21{1/4,\,-5/12}{1/3}{x} &\equal& { \left(1-\frac{5x}4 \right)^{1/4}}\;
\hpg{2}{1}{1/4,\,-1/12}{2/3}{-\frac{x^2(4x-5)^3}{(5x-4)^3}}, \\
\hpg21{-1/4,\,-7/12}{2/3}{x} &\equal& { \left(1-42x-7x^2 \right)^{1/4}}\;
\hpg{2}{1}{1/4,\,-1/12}{2/3}{-\frac{x\,(x^2-42x-7)^3}{(7x^2+42x-1)^3}}.
\end{eqnarray}
Or similarly, let $z=\varphi_{14}(x)$ be the degree 14 covering for the $(2/3,4/3,4/3)$
tetrahedral case. A hypergeometric identity is:
\begin{eqnarray}
\hpg{2}{1}{\!-1/2,-7/6}{-1/3}{x}=
\frac{(189x^4-378x^3+301x^2-112x+16)^{1/4}}{2}\;
\hpg{2}{1}{\!1/4,-1/12}{2/3}{\varphi_{14}(x)}.
\end{eqnarray}

A list of other transformations between algebraic Gauss hypergeometric equations is not long.
If the starting equation is not a standard (tetrahedral, octahedral or icosahedral) equation,
at least one of the local exponent differences is effectively non-restricted, so
we can only have special cases of classical transformations.
If one of the possible monodromy groups 
(i.e., icosahedral, octahedral, tetrahedral, dihedral or cyclic, including trivial) 
can be a subgroup of another, there is a transformation between two 
standard hypergeometric equations with those monodromy groups.
These transformations factor following the possible subgroup relations 
between the monodromy groups. The transformations that reduce
the monodromy group to a largest proper subgroup are special cases of transformations
considered in Sections \ref{clalgtr} through %, \ref{logarithms}, 
\ref{dihedrals}, except the transformation between standard icosahedral and
tetrahedral hypergeometric equations. The pull-back covering has degree 5:
\[
\varphi_5(x)=\frac{50(5\!+\!3\sqrt{-15})\,x\,(1024x\!-\!781\!-\!171\sqrt{-15})^3}
{(128x+7+33\sqrt{-15})^5}.
\]
Here is a corresponding hypergeometric identity:
\begin{equation} \label{icos2tetra}
\hpg{2}{1}{\!1/4,-1/12}{2/3}{\,x}= \left(
1\!+\!\frac{7\!-\!33\sqrt{-15}}{128}\,x\right)^{1/12} \hpg{2}{1}
{\!11/60,-1/60}{2/3}{\varphi_5(x)}.
\end{equation}
This transformation is derived in \cite[Section 5.1]{andkitaev} as well.
In general, if a standard hypergeometric equation is transformed to a
(not necessary standard) hypergeometric equation with smaller monodromy group,
that transformation must factor via the corresponding transformation between standard
equations and Klein's transformation preserving the smaller monodromy group.

\begin{remark} \label{trivexcept} \rm  
Hypergeometric equations with a finite monodromy group can be pull-backed 
to differential equations  with the trivial monodromy group. 
Then algebraic Gauss hypergeometric functions are transformed to rational functions 
(perhaps on a higher genus curve). The minimal transformation degree for 
these transformations is the order of the monodromy group, 
which is 12, 24, 60 for tetrahedral, octahedral, icosahedral equations, respectively.
If the transformed equation is hypergeometric, its degree is given by formula (\ref{trareas}). 

A priori, it seems possible that a pull-back to hypergeometric equation with
the trivial monodromy group can have all its 3 singularities outside the fibers above
$\{0,1,\infty\}\subset\PP_z^1$. This situation would be an exception to part 1 
of Lemma \ref{transeqv}: we would have a transformation between hypergeometric equations
without two-term identities between their hypergeometric solutions. 
A simple candidate for such a pull-back transformation could transform
the local exponent differences as $(1/2,1/2,1/2)\stackrel{\,10}{\longleftarrow}(2,2,2)$;
it would have 5 simple ramification points above each of the 3 points with the local exponent
difference $1/2$. 
It two hypergeometric equations with the trivial monodromy group can be transformed to each other,
they are related by a chain of transformations considered in Sections \ref{clalgtr} through 
\ref{alggtranfs}. In Klein's standard hypergeometric equations are involved,
there is a 
unique\footnote{Wrong. See \cite[Subsection 5.3]{dihedraltr} or the previous footnote %
for a counterexample.%

Let $E(a,b,c)$ denote a hypergeometric equation with the local exponent differences $a,b,c$.

The question of existence of pull-back transformations of hypergeometric functions %
that do not yield a two term-identity between their solutions is resolved positively %
in \cite[Remark 5.9]{dihedraltr}. For an example, consider a degree 12 composition %  
of transformations from $E(1/2,1/2,1/2)$ to $E(1,1,1)$, degree 4, %
and a general cubic transformation from $E(1,1,1)$ to $E(3,2,2)$. The transformation %
considered in Remark \ref{trivexcept} does not exist because $E(2,2,2)$ has logarithmic singularities %
rather than a trivial monodromy. }  %
(because of identification of Schwarz maps) transformation, 
which is a composition of considered transformations of standard equations reducing 
the projective monodromy group with Klein's transformation keeping the smallest monodromy group.
Otherwise we have a classical transformation. As we observed, all these transformations
allow two-term hypergeometric identities. (Particularly see the statement just before Remark 
\ref{analcont}). In the indicated composition of pull-back transformations, we can still have
a two-term identity if we keep a fractional local exponent difference at $z=0$ up till the last
transformation (possibly acting on a hypergeometric equation with the trivial monodromy group).
In particular, computations confirm that no pull-back covering for the transformation
$(1/2,1/2,1/2)\stackrel{\,10}{\longleftarrow}(2,2,2)$ exists.
\end{remark}

\section{Elliptic integrals} \label{ellints}

Here we consider algebraic transformations for solutions of 
hypergeometric equations with the local exponent differences $(1/2,1/4,1/4)$, or
$(1/2,1/3,1/6)$, or $(1/3,1/3,1/3)$. For each of these equations some of
hypergeometric solutions are trivial (i.e., constants or power functions), while
other solutions are incomplete elliptic integrals (up to a possible power
factor). Here are representative interesting solutions of hypergeometric equations
with the mentioned triples of local exponent differences: 
\begin{eqnarray} \label{elin4}
\hpg{2}{1}{1/2,\,1/4\,}{5/4}{\,z} & = & \frac{z^{-1/4}}{4}\, \int_0^z
t^{-3/4}\,(1-t)^{-1/2}\,dt \nonumber \\ \label{elin4a}
&=& \frac{z^{-1/4}}{2}\;%\int_0^{\sqrt{z}}
\int_{1/\!\sqrt{z}}^{\infty}\;\frac{dx}{\sqrt{x^3-x}}, \\
 \label{elin6} \hpg{2}{1}{1/2,\,1/6\,}{7/6}{\,z} &
= & \frac{z^{-1/6}}{6}\, \int_0^z t^{-5/6}\,(1-t)^{-1/2}\,dt \nonumber \\ \label{elin6a}
&=&\frac{z^{-1/6}}{2}%\int_0^{\sqrt[3]{z}}
\,\int_{1/\!\sqrt[3]{z}}^{\infty}\;\frac{dx}{\sqrt{x^3-1}},\\
\label{elin3} \hpg{2}{1}{1/3,\,2/3\,}{4/3}{\,z} & = & \frac{z^{-1/3}}{3}\,
\int_0^z t^{-2/3}\,(1-t)^{-2/3}\,dt \nonumber  \\ \label{elin3a}
&=& z^{-1/3}\,\int_{1/\!\sqrt[3]{z}}^\infty\;\frac{dx}{(x^3-1)^{2/3}}.
\end{eqnarray}
Here we substituted $t=x^{-2}$ or $t=x^{-3}$ into the immediate
integral expressions. % in (\ref{elin4}), (\ref{elin6}), (\ref{elin3}).
As we see, the three hypergeometric functions can be transformed to integrals 
of holomorphic forms on the genus 1 curves 
\begin{eqnarray} \label{eq:ecurves}
y^2=x^3-x, \qquad y^2=x^3-1,\qquad x^3+y^3=1,
\end{eqnarray}
respectively. In fact, the integrand functions define algebraic curves isomorphic respectively 
to these three cubic curves; see \cite{telliptici} for details. 
Let $E_1$, $E_2$, $E_3$ be three curves defined in (\ref{eq:ecurves}),
respectively. We consider them as elliptic curves (with the classical group structure) 
by fixing the point at infinity for $E_1$ and $E_2$,  or the infinite point $(1\!:\!-1\!:\!0)$ for $E_3$,
as the neutral element of the group structure. The elliptic curves $E_2$ and $E_3$ are isomorphic.
Non-trivial solutions of hypergeometric equations with the local exponent differences
$(2/3,1/6,1/6)$ are genus 2 hyperelliptic integrals. For example,
\begin{eqnarray*} \label{hpelin}
\hpg{2}{1}{1/3,\,1/6}{7/6}{z} & = & 
\frac{z^{-1/6}}{6}\, \int_0^zt^{-5/6}\,(1-t)^{-1/3}\,dt \nonumber \\
& = & \frac{z^{-1/6}}{2^{1/3}}\,
\int_{\theta(z)}^{\infty} \frac{X\,dX}{\sqrt{X^6+1}},
\qquad\mbox{where}\quad
\theta(z)=\frac{(1-z)^{1/3}}{2^{1/3}\,z^{1/6}}.
\end{eqnarray*}
Here the substitution is $t\to \left(\sqrt{X^6+1}-X^3\right)^2$.

Formula (\ref{trareas}) gives no restriction on the degree $d$ of pull-back transformations of the hypergeometric equations under consideration. But (\ref{trareas}) requires that the transformed hypergeometric equation must have %positive
local exponent differences $e_0,e_1,e_\infty$ such that
\mbox{$e_0+e_1+e_\infty=1$}. In particular, all 3 singularities of the
transformed equation are relevant singularities, 
and the pull-back covering branches only above 3 points. 
Possible branching patterns are presented in Table \ref{elltab}.
Multiplicative terms in the last column give the branching order (as the second multiplicant)
and the number of points with that branching order in the same fiber (as the first multiplicant).
\begin{table}
\begin{center} \begin{tabular}{|c|c|c|c|}
\hline \multicolumn{2}{|c|}{Local exponent differences} & Degree &
Branching above \\ \cline{1-2} %$\!(1/k_1,\,1/k_2,\,1/k_3)\!$
below & above & $d$ & the regular singular points \\
 \hline
$(1/2,\,1/4,\,1/4)$ & $(1/2,\,1/4,\,1/4)$ & $4n$ & $2n\!*2=n\!*4=(n\!-\!1)\!*4+2+1+1$ \\
$(1/2,\,1/4,\,1/4)$ & $(1/2,\,1/4,\,1/4)$ & $4n\!+\!1$ & $2n\!*2+1=n\!*4+1=n\!*4+1$ \\
$(1/2,\,1/4,\,1/4)$ & $(1/2,\,1/4,\,1/4)$ & $4n\!+\!2$ & $(2n\!+\!1)\!*2=n\!*4+2=n\!*4+1+1$ \\
$(1/2,\,1/3,\,1/6)$ & $(1/2,\,1/3,\,1/6)$ & $6n$ & $3n\!*2=2n\!*3=(n\!-\!1)\!*6+3+2+1$ \\
$(1/2,\,1/3,\,1/6)$ & $(1/2,\,1/3,\,1/6)$ & $6n\!+\!1$ & $3n\!*2+1=2n\!*3+1=n\!*6+1$ \\
$(1/2,\,1/3,\,1/6)$ & $(1/2,\,1/3,\,1/6)$ & $6n\!+\!3$ &
$(3n\!+\!1)\!*2\!+\!1=(2n\!+\!1)\!*3=n\!*6\!+\!2\!+\!1$ \\
$(1/2,\,1/3,\,1/6)$ & $(1/2,\,1/3,\,1/6)$ & $6n\!+\!4$ &
$(3n\!+\!2)\!*2=(2n\!+\!1)\!*3\!+\!1=n\!*6\!+\!3\!+\!1$ \\
$(1/2,\,1/3,\,1/6)$ & $(1/3,\,1/3,\,1/3)$ & $6n$ & $3n\!*2=2n\!*3=(n\!-\!1)\!*6+2+2+2$ \\
$(1/2,\,1/3,\,1/6)$ & $(1/3,\,1/3,\,1/3)$ & $6n$ & $3n\!*2=(2n\!-\!1)\!*3+1+1+1=n\!*6$ \\
$(1/2,\,1/3,\,1/6)$ & $(1/3,\,1/3,\,1/3)$ & $6n\!+\!2$ & $(3n\!+\!1)\!*2=2n\!*3+1+1=n\!*6+2$ \\
$(1/2,\,1/3,\,1/6)$ & $(1/3,\,1/3,\,1/3)$ & $6n\!+\!4$ &
$(3n\!+\!2)\!*2=(2n\!+\!1)\!*3\!+\!1=n\!*6\!+\!2\!+\!2$ \\
$(1/2,\,1/3,\,1/6)$ & $(2/3,\,1/6,\,1/6)$ & $6n$ & $3n\!*2=2n\!*3=(n\!-\!1)\!*6+4+1+1$ \\
$(1/2,\,1/3,\,1/6)$ & $(2/3,\,1/6,\,1/6)$ & $6n\!+\!2$ & $(3n\!+\!1)\!*2=2n\!*3+2=n\!*6+1+1$ \\
$(1/3,\,1/3,\,1/3)$ & $(1/3,\,1/3,\,1/3)$ &
$3n$ & $n\!*3=n\!*3=(n\!-\!1)\!*3+1+1+1$ \\
$(1/3,\,1/3,\,1/3)$ & $(1/3,\,1/3,\,1/3)$ & $3n\!+\!1$ &
$n\!*3\!+\!1=n\!*3\!+\!1=n\!*3\!+\!1$ \\ \hline
\end{tabular} \end{center}
\caption{Transformations of hypergeometric elliptic integrals}
\label{elltab}
\end{table}

Coverings with the branching patterns of Table \ref{elltab} give rise to morphisms
between the corresponding (hyper)elliptic curves. For example,
a pull-back transformation $(1/2,1/3,1/6)\leftarrow(1/3,1/3,1/3)$ of degree $6n+2$
implies the polynomial identity
\begin{equation} \label{pqrcheck2}
R_{3n+1}(z)^2=(z-1)Q_{2n}(z)^3-z^2\,P_n(z)^6 
\end{equation}
for some polynomials $P_n(z)$, $Q_{2n}(z)$, $R_{3n+1}(z)$ 
of degree $n$, $2n$, $3n+1$, respectively. This gives the following
morphism from $E_3$ to $E_2$:
\begin{equation}
(x,y)\mapsto \left(
\frac{x\,y\,Q_{2n}(x^{-3})}{P_n(x^{-3})^2},\;
\frac{x^3\,R_{3n+1}(x^{-3})}{P_n(x^{-3})^3} \right).
\end{equation}
Conversely, a morphism (or endomorphism) between the elliptic curves relates
the holomorphic differentials up to a constant multiple, and gives rise to a 
transformation of their integrals. If the morphism fixes the upper integration bound $\infty$ 
in (\ref{elin4a}), (\ref{elin6a}), (\ref{elin3a}), we get a transformation of the
hypergeometric functions as well. This correspondence is investigated thoroughly in \cite{telliptici}.
Here we demonstrate it on several examples.

Most interesting are pull-back transformations of the three mentioned
hypergeometric equations to themselves. %As we demonstrate here,
These transformations correspond to isogeny endomorphisms 
of the elliptic curves $E_1$, $E_2$ or $E_3$. 
%This observation was immediately suggested by Beukers \cite{beukrspriv}. 
The ring of isogeny endomorphisms for $E_1$ is isomorphic to the ring $\ZZ[i]$ of Gaussian integers 
\cite{silverman1}. The ring of isogeny endomorphisms for $E_2$ or $E_3$  is isomorphic to the ring 
$\ZZ[\omega]$, where $\omega$ is a primitive cubic root of unity as in (\ref{cubic3}).
Composition of the isogenies corresponds to multiplication in the mentioned rings 
of algebraic integers. Roots of unity in both $\ZZ[i]$ and $\ZZ[\omega]$ correspond to trivial
transformations of hypergeometric equations. The degree of a pull-back transformation
induced by an endomorphism is equal to the $\CC$-norm of the corresponding algebraic integer.
In particular, there may be none or several transformations for a fixed degree and branching pattern
from Table \ref{elltab}, depending on how many algebraic integers exist with that norm.

The transformations from the local exponent differences $(1/2,1/3,1/6)$ to 
the local exponent differences $(1/3,1/3,1/3)$ or $(2/3,1/6,1/6)$
are compositions of a classical quadratic transformation and the mentioned 
endomorphisms of elliptic curves; see \cite[Section 5]{telliptici}.
In particular, there are actually no transformations 
$(1/2,1/3,1/6)\leftarrow(1/3,1/3,1/3)$ of degree $6n+4$,
even if indicated in Table \ref{elltab}, because there are no transformations
$(1/2,1/3,1/6)\leftarrow(1/2,1/3,1/6)$ of degree $3n+2$.

Now we consider explicitly pull-back transformations coming from the endomorphisms of $E_1$.
If $(x,y)\mapsto(\psi_x,\psi_y)$ is an isogeny endomorphism of $E_1$, then the substitution
$x\mapsto\psi_x(x,\sqrt{x^3-x})$ into (\ref{elin4a}) gives an integral of a
holomorphic differential form again. Since the linear space of holomorphic
differentials on $E_1$ is one-dimensional, the transformed differential form
must be proportional to $dx/\sqrt{x^3-x}$. The upper integration bound does
not change. Transformation of the lower integration bound gives the
transformation $z\mapsto\psi_x(1/\sqrt{z})^{-2}$ of the hypergeometric
function into itself, up to a radical factor. Using induction and the addition law on $E_1$,
one can prove \cite[Theorem 2.1]{telliptici} that $\psi_x(1/\sqrt{z})^{-2}$ is a rational function
for any isogeny endomorphism, and that its degree is equal to the degree of the isogeny.
Conversely \cite[Theorem 2.1]{telliptici}, analysis of the first three branching patterns 
in Table \ref{elltab} shows  that any pull-back transformation 
$(1/2,1/4,1/4)\stackrel{d}{\longleftarrow}(1/2,1/4,1/4)$ is 
induced by an endomorphism of $E_1$.

As mentioned, the ring of isogeny endomorphisms of $E_1$ is isomorphic to the ring $\ZZ[i]$ of
Gaussian integers. We identify $i\in\ZZ[i]$ with the isogeny $(x,y)\mapsto (-x,iy)$. 
Addition of isogenies is equivalent to the chord-and-tangent addition law on $E_1$. 
Here are a few examples of isogenies on $E_1$:
\begin{eqnarray*}
(x,y) \mapsto \left(\frac{x^2-1}{2i\,x},
\;\frac{y\;(x^2+1)}{2(i-1)\,x^2}\right),\qquad
(x,y)\mapsto\left(\frac{(x^2+1)^2}{4\,x\,(x^2-1)},
\frac{(x^2+1)(x^4-6x^2+1)}{8\,x\,y\,(x^2-1)}\right),\\
(x,y)\mapsto\left( \frac{x(x^2\!-\!1\!-\!2i)^2}{((1\!+\!2i)x^2\!-\!1)^2},
\frac{y(x^4\!+\!(2\!+\!8i)x^2\!+\!1)(x^2\!-\!1\!-\!2i)}
{((1\!+\!2i)x^2\!-\!1)^3}\right).
\end{eqnarray*}
They correspond to the Gaussian integers $1+i$, $2$, $1+2i$, respectively.
Here below are the induced algebraic transformations of Gauss hypergeometric functions:
\begin{eqnarray}
\hpg{2}{1}{1/2,\,1/4}{5/4}{\,z} & = & \frac{1}{\sqrt{1-z}}\;
\hpg{2}{1}{1/2,\,1/4}{5/4}{-\frac{4\,z}{(z-1)^2}}.\\
\hpg{2}{1}{1/2,\,1/4}{5/4}{\,z} & = & \frac{\sqrt{1-z}}{1+z}\;
\hpg{2}{1}{1/2,\,1/4}{5/4}{\frac{16\,z\,(z-1)^2}{(z+1)^4}}.\\ \label{elltr5}
\hpg{2}{1}{1/2,\,1/4}{5/4}{\,z}&=&\frac{1-z/(1\!+\!2i)}{1-(1\!+\!2i)z}\;
\hpg{2}{1}{1/2,\,1/4}{5/4}{\frac{z\,(z-1-2i)^4}{\big((1\!+\!2i)z-1\big)^4}}.
\end{eqnarray}
The first two identities are special cases of classical transformations. The radical
factors on a right-hand side are equal to
$\left[\psi_x(1/\sqrt{z})^{-2}\right]^{1/4}z^{-1/4}$ times a constant 
(deducible from the value of hypergeometric series at $z=0$).
The transformations $(1/2,1/4,1/4)\stackrel{d}{\longleftarrow}(1/2,1/4,1/4)$
form a semi-group under composition, isomorphic to the multiplicative semi-group 
$\ZZ[i]^{\,\star}/(\pm 1,\pm i)$. 
% the transformation group is isomorphic to $\ZZ[i]^{\,\star}/(\pm 1,\pm i)$,
%since the units of $\ZZ[i]$ give trivial transformations. 
The degree of these transformations is equal to the norm $p^2+q^2$ of
a corresponding Gaussian integer $p+qi$. In particular, 
there are no transformations of degree $21$ although Table \ref{elltab} allows it, 
because there are no Gaussian integers with this norm. On the other hand, 
there are several different transformations of degree 25, corresponding to $3\pm4i$ or $5$. 
One of them is the composition of (\ref{elltr5}) with itself, the other is the
composition of (\ref{elltr5}) with the complex conjugate of itself. In fact,
algebraic transformations related by the complex conjugation are not related by
fractional-linear transformations in general. The addition law on $E_1$ can be translated 
into ``addition" of the polynomial triples determining the branching points (of order 2 or 4)
of explicit pull-back coverings for the first three branching patterns in Table \ref{elltab};
see \cite[Section 2]{telliptici}.

Likewise, an isogeny endomorphism on $E_2$ 
transforms the holomorphic differential form in (\ref{elin6a}) into a
scalar multiple of itself, and the upper integration bound does not change.
The lower integration bound changes as $z\mapsto\psi_x(z^{-1/3})^{-3}$.
By induction and the addition law on $E_2$, this is a rational function
determining a desired pull-back covering, and its degree is equal to the degree of the isogeny
\cite[Section 3]{telliptici}. Conversely, analysis of respective cases 
of Table \ref{elltab} shows  that any pull-back transformation 
$(1/2,1/3,1/6)\stackrel{d}{\longleftarrow}(1/2,1/3,1/6)$ is 
induced by an endomorphism of $E_2$.
These transformations form a semi-group under composition, 
isomorphic to $\ZZ[\omega]^{\,\star}/\left(\pm 1,\pm\omega,\pm\omega^{-1}\right)$. 
We identify the cubic root $\omega$ with the isogeny $(x,y)\mapsto (\omega x,y)$.
Here are examples of explicit transformations corresponding to the algebraic
integers $1-\omega$, $3$, $3+\omega$ of $\ZZ[\omega]$:
\begin{eqnarray}
\hpg{2}{1}{\!1/2,\,1/6}{7/6}{\,z} & \equal & \frac{1}{\sqrt{1-4z}}\;
\hpg{2}{1}{\!1/2,\,1/6}{7/6}{\frac{27\,z}{(4z-1)^3}},\\
\hpg{2}{1}{\!1/2,\,1/6}{7/6}{\,z} & \equal &
\frac{1-4z}{\sqrt{1\!+\!96z\!+\!48z^2\!-\!64z^3}}\,\hpg{2}{1}{\!1/2,1/6}{7/6}
{\frac{-729\,z\,(4z-1)^6}{(64z^3\!-\!48z^2\!-\!96z\!-\!1)^3}},\\
\hpg{2}{1}{\!1/2,1/6}{7/6}{\,z} &\equal&\frac{1-4z/(3\omega\!+\!1)}
{\sqrt{1\!-\!(44\!+\!48\omega)z\!+\!(48\omega\!+\!16)z^2}}\nonumber\\
&&\times\,\hpg{2}{1}{\!1/2,1/6}{7/6}{\frac{z\;(4z-3\omega\!-\!1)^6}
{((48\omega\!+\!16)z^2\!-\!(44\!+\!48\omega)z\!+\!1)^3}}.
\end{eqnarray}

Similarly, transformations $(1/3,1/3,1/3)\stackrel{d}{\longleftarrow}(1/3,1/3,1/3)$ 
correspond to the isogeny endomorphisms on $E_3$. Recall that this elliptic curve
is isomorphic to $E_2$. With the chosen addition law on $E_3$ and identification 
of a hypergeometric solution as the integral in (\ref{elin3a}), the isogeny of
multiplication by $-1\in\ZZ[\omega]$ corresponds to Euler's transformation
(\ref{flinear1}). Transformations of the hypergeometric function for (\ref{elin3a})
into itself form a semi-group (under composition) isomorphic to 
$\ZZ[\omega]^*/\left(1,\omega,\omega^{-1}\right)$. We identify the cubic root
$\omega$ with the isogeny $(x,y)\mapsto(\omega^{-1}x,\omega^{-1}y)$.
Here are explicit transformations corresponding to $2$,  % $1-\omega$, 
$3$,  $3+\omega\in\ZZ[\omega]$:
\begin{eqnarray}
%\hpg{2}{1}{1/3,2/3}{4/3}{\,z} & \equal &
%\frac{(1-z)^{1/3}}{1+\omega^2z}\;\hpg{2}{1}{1/3,\,2/3}{4/3}
%{\frac{3\,(2w+1)\,z\,(z-1)}{(z+\omega)^3}}.\\
\hpg{2}{1}{1/3,2/3}{4/3}{\,z} & \equal &
\frac{1-z/2}{1-2z}\;\hpg{2}{1}{1/3,\,2/3}{4/3}{\frac{z\,(z-2)^3}{(1-2z)^3}},\\
\hpg{2}{1}{1/3,2/3}{4/3}{\,z} & \equal &
\frac{(1\!-\!z\!+\!z^2)\,(1\!-\!z)^{1/3}}
{1\!+\!3z\!-\!6z^2\!+\!z^3}\;\hpg{2}{1}{1/3,\,2/3}{4/3}
{\frac{27\,z\,(z-1)\,(z^2-z+1)^3}{(z^3-6z^2+3z+1)^3}}.\\
\hpg{2}{1}{\!1/3,2/3}{4/3}{\,z}&\equal&
\frac{1-z-z^2/(3\omega\!+\!2)}{1\!+\!(3\omega\!+\!2)z\!-\!(3\omega\!+\!2)z^2}
\,\hpg{2}{1}{\!1/3,2/3}{4/3}{\frac{z\;(z^2+(3\omega\!+\!2)z-3\omega\!-\!2)^3}
{(1+(3\omega\!+\!2)z-(3\omega\!+\!2)z^2)^3}}.
\end{eqnarray}
An explicit transformation corresponding to $1-\omega$ can be obtained from (\ref{cubic3})
with $a=1$.

As mentioned, transformations from the local exponent differences $(1/2,1/3,1/6)$ to 
the local exponent differences $(1/3,1/3,1/3)$ or $(2/3,1/6,1/6)$
are compositions of a classical quadratic transformation and the mentioned 
endomorphisms of elliptic curves.
Correspondingly, the morphisms from $E_3$ or the hyperelliptic curve $Y^2=X^6+1$
to $E_2$, that leave the infinite points at infinity, factor via isogeny
endomorphisms of $E_2$ and the straightforward morphisms
$(x,y)\mapsto (2^{2/3}xy,i-2ix^3)$ or $(X,Y)\mapsto (-X^2,iY)$,
%\begin{eqnarray*}
%(x,y)\mapsto (2^{2/3}xy,i-2ix^3) \qquad\mbox{or}\qquad
%(X,Y)\mapsto (-X^2,iY),
%\end{eqnarray*}
respectively.

\section{Hyperbolic hypergeometric functions}
\label{otherat}

Transformations of hyperbolic hypergeometric functions are extensively
studied in \cite{thyperbolic}. There are 9 non-classical transformations in this case,
of degree 6, 8, 9, 10, 12, 18 or 24.

Without loss of generality, we may assume $k_1\le k_2\le k_3\le d$. 
Inequality (\ref{nsingpoints2}) together with
$1/k_1+1/k_2+1/k_3<1$ already implies finitely many possibilities for the
tuple $(k_1,k_2,k_3,d)$. Indeed, inequality (\ref{nsingpoints2}) gives a
bound for $d$ when $k_1,k_2,k_3$ are fixed; then $k_3\le d$ gives a bound
for $k_3$ when $k_1,k_2$ are fixed, etc. But stronger inequalities
and conditions follow from \cite[Lemma 2.2]{thyperbolic}. First of
all, the transformed equation must have precisely 3 singular points, and the
covering $\varphi:\PP^1_z\to\PP_x^1$ branches only above the set
$\{0,1,\infty\}\subset\PP_z^1$. Then we consequently derive:
\begin{eqnarray*}
d-\left\lfloor \frac{d}{k_1}\right\rfloor-\left\lfloor \frac{d}{k_2}
\right\rfloor-\left\lfloor\frac{d}{k_3}\right\rfloor=1,&\displaystyle\quad\,
d\left(1-\frac{1}{k_1}-\frac{1}{k_2}-\frac{1}{k_3}\right)\le 1-\frac{3}{k_3},\\
\left(1-\frac1{k_1}-\frac1{k_2}\right)k_3^2-2k_3+3\le 0,
&\displaystyle\qquad\quad\!\frac23\le \frac1{k_1}+\frac1{k_2}<1.
\end{eqnarray*}
With these stronger formulas we get a moderate list of possibilities after
Step 2 of our classification scheme. The list of possible branching patterns
after Step 3 is presented by the first three columns of Table \ref{figtab}.
The branching patterns are determined by the two triples of local exponent
differences and the principle that each fiber of
$\{0,1,\infty\}\subset\PP_z^1$ contains maximal possible number of
non-singular points. For each branching pattern there is at most one
covering. The coverings were computed by the algorithm in \cite[Section
3]{thyperbolic}; they are characterized in the fourth column of Table
\ref{figtab}. The last column indicates existence of Coxeter decompositions
described at the end of Section \ref{clscheme}. The three cases which admit
a Coxeter decomposition are implied in \cite{hodgkins2} and \cite{beukers}.
\begin{table}
\begin{center} \begin{tabular}{|c|c|c|c|c|}
\hline \multicolumn{2}{|c|}{Local exponent differences} & Degree & Covering
& Coxeter \\ \cline{1-2}
$(1/k,\,1/\ell,\,1/m)$ & above & $d$ & composition & decomposition \\
 \hline
$(1/2,\,1/3,\,1/7)$ & $(1/3,\,1/3,\,1/7)$ & 8 & indecomposable &  no \\
$(1/2,\,1/3,\,1/7)$ & $(1/2,\,1/7,\,1/7)$ & 9 & indecomposable & no \\
$(1/2,\,1/3,\,1/7)$ & $(1/3,\,1/7,\,2/7)$ & 10 & indecomposable & yes \\
$(1/2,\,1/3,\,1/7)$ & $(1/7,\,1/7,\,3/7)$&12&\multicolumn{2}{c|}{no covering}\\
$(1/2,\,1/3,\,1/7)$ & $(1/7,\,2/7,\,2/7)$&12&\multicolumn{2}{c|}{no covering}\\
$(1/2,\,1/3,\,1/7)$ & $(1/3,\,1/7,\,1/7)$&16&\multicolumn{2}{c|}{no covering}\\
$(1/2,\,1/3,\,1/7)$ & $(1/7,\,1/7,\,2/7)$ & 18 & $2\times 9$ & no \\
$(1/2,\,1/3,\,1/7)$ & $(1/7,\,1/7,\,1/7)$ & 24 & $3\times 8$ & yes \\
$(1/2,\,1/3,\,1/8)$ & $(1/3,\,1/8,\,1/8)$ & 10 & indecomposable & no \\
$(1/2,\,1/3,\,1/8)$ & $(1/4,\,1/8,\,1/8)$ & 12 & $2\times 2\times 3$ & yes \\
$(1/2,\,1/3,\,1/9)$ & $(1/9,\,1/9,\,1/9)$ & 12 & $3\times 4$ & no \\
$(1/2,\,1/4,\,1/5)$ & $(1/4,\,1/4,\,1/5)$ & 6 & indecomposable & no \\
$(1/2,\,1/4,\,1/5)$ & $(1/5,\,1/5,\,1/5)$ & 8 &\multicolumn{2}{c|}{no
covering}\\ \hline
\end{tabular} \end{center}
\caption{Transformations of hyperbolic hypergeometric functions}
\label{figtab}
\end{table}

Here we give rational functions defining the indecomposable pull-back
transformations, and examples of corresponding algebraic transformations of
Gauss hypergeometric functions.
\begin{itemize}
\item $(1/2,\,1/3,\,1/7)\stackrel{8}{\longleftarrow}(1/3,\,1/3,\,1/7)$. This
transformation was independently computed numerically in \cite{andkitaev}, 
and later fully presented in \cite{kitaevdb}. Let $\omega$
denote a primitive cubic root of unity as in (\ref{cubic3}). The covering
and an algebraic transformation are:
\begin{eqnarray}
\varphi_8(x)&\equal&
\frac{x\,(x-1)\,\left(27x^2-(723\!+\!1392\omega)x-496\!+\!696\omega\right)^3}
{64\,\big((6\omega+3)x-8-3\omega\big)^7},\\
\hpg{2}{1}{2/21,\,5/21}{2/3}{\,x}&\equal&
{\textstyle\left(1-\frac{33+39\omega}{49}\,x\right)}^{-1/12}
\;\hpg{2}{1}{1/84,\,13/84}{2/3}{\varphi_8(x)}.
\end{eqnarray}
Note that the conjugation $\omega=-1-\omega$ acts in the same way as a
composition with fractional-linear transformation interchanging the points
$x=0$ and $x=1$. This confirms uniqueness of the covering.%
\item $(1/2,\,1/3,\,1/7)\stackrel{9}{\longleftarrow}(1/2,\,1/7,\,1/7)$. Let
$\xi$ denote an algebraic number satisfying $\xi^2+\xi+2=0$. The covering
and an algebraic transformation are:
\begin{eqnarray}
\varphi_9(x)&\equal&\frac{27\,x\,(x-1)\,(49x-31-13\xi)^7}
{49\,(7203x^3+(9947\xi-5831)x^2-(9947\xi+2009)x+275-87\xi)^3},\\
\hpg{2}{1}{3/28,\,17/28}{6/7}{x} &\equal&  \textstyle{\left(
1+\frac{7(10-29\xi)}{8}\,x-\frac{343(50-29\xi)}{512}\,x^2
+\frac{1029(362+87\xi)}{16384}\,x^3\right)}^{-1/28}\nonumber\\
&& \times\,\hpg{2}{1}{1/84,\,29/84}{6/7}{\varphi_9(x)}.
\end{eqnarray}
\item $(1/2,\,1/3,\,1/7)\stackrel{10}{\longleftarrow}(1/3,\,1/7,\,2/7)$. 
This transformation was independently computed in \cite{kitaevdb} as well.
The covering and an algebraic transformation are:
\begin{eqnarray}
\varphi_{10}(x)&\equal&-\frac{x^2\,(x-1)\,(49x-81)^7}{4\,(16807x^3-9261x^2-13851x+6561)^3},\\
\hpg{2}{1}{\!5/42,\,19/42}{5/7}{x} & \!\equal\! &
{\textstyle\left(1\!-\!\frac{19}9x\!-\!\frac{343}{243}x^2\!+\!\frac{16807}{6561}x^3\right)}^{-1/28}
\hpg{2}{1}{1/84,\,29/84}{6/7}{\varphi_{10}(x)}.
\end{eqnarray}
\item $(1/2,\,1/3,\,1/8)\stackrel{10}{\longleftarrow}(1/3,\,1/8,\,1/8)$. Let
$\beta$ denote an algebraic number satisfying $\beta^2+2=0$. The covering
and an algebraic transformation are:
\begin{eqnarray}
\widetilde{\varphi}_{10}(x)&\equal&\frac{4\,x\,(x-1)\,(8\beta x+7-4\beta)^8}
{(2048\beta x^3-3072\beta x^2-3264x^2+912\beta x+3264x+56\beta-17)^3},\\
\hpg{2}{1}{\!5/24,13/24}{7/8}{x} &\equal &
{\textstyle\left(1+\frac{16(4-17\beta)}9x-\frac{64(167-136\beta)}{243}x^2
+\frac{2048(112-17\beta)}{6561}x^3\right)}^{-1/16}\nonumber\\
&&\times\,\hpg{2}{1}{1/48,\,17/48}{7/8}{\widetilde{\varphi}_{10}(x)}.
\end{eqnarray}
\item $(1/2,\,1/4,\,1/5)\stackrel{6}{\longleftarrow}(1/4,\,1/4,\,1/5)$. The
covering and an algebraic transformation are:
\begin{eqnarray}
\varphi_6(x)&\equal&\frac{4i\,x\,(x-1)\,(4x-2-11i)^4}{(8x-4+3i)^5},\\
\!\hpg{2}{1}{\!3/20,\,7/20}{3/4}{x} &\equal&
{\textstyle\left(1-\frac{8(4+3i)}{25}x\right)}^{-1/8}
\hpg{2}{1}{1/40,\,9/40}{3/4}{\varphi_6(x)}.
\end{eqnarray}
\end{itemize}
The composite transformations can be schematically represented similarly as
in Section \ref{clalgtr}:
\begin{eqnarray*}
&(1/2,\,1/3,\,1/7)\stackrel{9}{\longleftarrow}(1/2,\,1/7,\,1/7)
\stackrel{2}{\longleftarrow}(1/7,\,1/7,\,2/7),\\
&(1/2,\,1/3,\,1/7)\stackrel{8}{\longleftarrow}(1/3,\,1/3,\,1/7)
\stackrel{3}{\longleftarrow}(1/7,\,1/7,\,1/7),\\
&(1/2,1/3,1/8)\stackrel{3}{\longleftarrow}(1/2,1/4,1/8)\stackrel{2}
{\longleftarrow}(1/2,1/8,1/8)\stackrel{2}{\longleftarrow}(1/4,1/8,1/8),\\
&(1/2,\,1/3,\,1/9)\stackrel{4}{\longleftarrow}(1/3,\,1/3,\,1/9)
\stackrel{3}{\longleftarrow}(1/9,\,1/9,\,1/9).
\end{eqnarray*}
Note that the transformation of degree 24 admits a Coxeter decomposition,
although it is a composition of two transformations without a Coxeter
decomposition. Here is an explicit algebraic transformation of degree 24:
\begin{eqnarray}
\hpg{2}{1}{\!2/7,\,3/7}{6/7}{x} &\equal &
\left(1-235x+1430x^2-1695x^3+270x^4+229x^5+x^6\right)^{-1/28}\nonumber\\
&&\times(1-x+x^2)^{-1/28}\;\hpg{2}{1}{1/84,\,29/84}{6/7}{\varphi_{24}(x)},
\end{eqnarray}
where $\displaystyle\varphi_{24}(x)=\frac{1728\,x\,(x-1)\,(x^3-8x^2+5x+1)^7}
{(x^2-x+1)^3(x^6+229x^5+270x^4-1695x^3+1430x^2-235x+1)^3}$.

\vspace{24pt}

\noindent{\bf Acknowledgements.} The author would like to thank Robert S.
Maier, Frits Beukers, Anton H.M.~Levelt and Masaaki Yoshida for useful
references and suggestions.

\bibliographystyle{alpha}
\bibliography{../hypergeometric}

\end{document}